\documentclass{amsart}
\usepackage{amsmath}
\usepackage{amssymb}
\usepackage{graphicx}
\usepackage{color}
\usepackage[all]{xy}
\usepackage{amsthm}
\theoremstyle{plain}

\newtheorem{thm}{Theorem}[section]
\newtheorem*{thm*}{Theorem}
\newtheorem{lem}[thm]{Lemma}
\newtheorem*{lem*}{Lemma}
\newtheorem{cor}[thm]{Corollary}
\newtheorem*{cor*}{Corollary}

\newtheorem*{cla*}{Claim}
\newtheorem{pro}[thm]{Proposition}
\newtheorem*{pro*}{Proposition}

\newtheorem*{fac*}{Fact}

\newtheorem*{que*}{Question}
\newtheorem{prob}[thm]{Problem}
\newtheorem*{prob*}{Problem}
\newtheorem{con}[thm]{Conjecture}
\newtheorem*{con*}{Conjecture}

\newtheorem*{rem*}{Remark}

\newtheorem*{rems*}{Remarks}

\newtheorem*{defn*}{Definition}

\newcommand{\R}{\ensuremath{\mathbb{R}}}
\newcommand{\C}{\ensuremath{\mathbb{C}}}

\newcommand{\Z}{\ensuremath{\mathbb{Z}}}
\newcommand{\N}{\ensuremath{\mathbb{N}}}
\newcommand{\Q}{\ensuremath{\mathbb{Q}}}

\newcommand{\pslc}{\mathrm{PSL}_2\C}

\newcommand{\tr}{\mathrm{tr}}
\newcommand{\im}{\mathrm{Im}}
\newcommand{\re}{\mathrm{Re}}

\begin{document}

\title{J\o{}rgensen Number and Arithmeticity}
\author{Jason Callahan}
\maketitle

\begin{abstract}
The J\o{}rgensen number of a rank-two non-elementary Kleinian group $\Gamma$ is
\[ J(\Gamma) = \inf\{|\tr^2 X - 4| + |\tr [X, Y] - 2| : \langle X, Y \rangle = \Gamma \}. \]
J\o{}rgensen's Inequality guarantees $J(\Gamma) \geq 1$, and $\Gamma$ is a J\o{}rgensen group if $J(\Gamma) = 1$.  This paper shows that the only torsion-free J\o{}rgensen group is the figure-eight knot group, identifies all non-cocompact arithmetic J\o{}rgensen groups, and establishes a characterization of cocompact arithmetic J\o{}rgensen groups. The paper concludes with computations of $J(\Gamma)$ for several non-cocompact Kleinian groups including some two-bridge knot and link groups.
\end{abstract}

\section{Introduction}

In \cite{Jorgensen76}, T.\ J\o{}rgensen establishes the following well known necessary condition for two elements of $\pslc$ to generate a non-elementary discrete group.

\newtheorem*{ji}{J\o{}rgensen's Inequality}
\begin{ji}
If $\langle X, Y \rangle$ is a non-elementary Kleinian group, then 
\[ |\tr^2 X - 4| + |\tr [X, Y] - 2| \geq 1. \]
\end{ji}

Accordingly, the \emph{J\o{}rgensen number} of an ordered pair of elements in $\pslc$ is
\[ J(X, Y) = |\tr^2 X - 4| + |\tr [X, Y] - 2|, \]
and the \emph{J\o{}rgensen number} of a rank-two non-elementary Kleinian group $\Gamma$ is
\[ J(\Gamma) = \inf\{J(X, Y) : \langle X, Y \rangle = \Gamma \}. \]
J\o{}rgensen's Inequality guarantees $J(\Gamma) \geq 1$, and if $J(\Gamma) = 1$, then $\Gamma$ is a \emph{J\o{}rgensen group} and has been the subject of much study.  Among the first such results is the following combination of Theorems 1 and 2 in \cite{JorgensenKiikka75}.

\begin{thm} \label{kiikka}
If $\langle X,Y \rangle$ is a non-elementary Kleinian group with $J(X,Y) = 1$, then the following two statements hold.
\begin{enumerate}
 \item $\langle X, YXY^{-1} \rangle$ is also non-elementary and $J(X, YXY^{-1}) = 1$.
 \item $X$ is parabolic, or $X$ is elliptic of order at least seven and $\tr XYXY^{-1} = 1$.
\end{enumerate}
\end{thm}

A J\o{}rgensen group $\Gamma = \langle X,Y \rangle$ with $J(X,Y) = 1$ is then \emph{of parabolic} or \emph{elliptic type} according to whether $X$ is parabolic or elliptic as a consequence of Theorem \ref{kiikka}.

The following problem has attracted much attention (e.g., \cite{JorgensenKiikka75}, \cite{GehringMartin89}, \cite{SatoYamada93}, \cite{Sato00}, \cite{LiOichiSato04}, \cite{LiOichiSato05a}, \cite{LiOichiSato05b}, and \cite{Gonzalez-AcunaRamirez07}).

\begin{prob} \label{prb}
Identify all J\o{}rgensen groups.
\end{prob}

It is observed in \cite{JorgensenKiikka75} that there are uncountably many non-conjugate J\o{}rgensen groups in general, so most work on Problem \ref{prb} has entailed restriction to more tractable cases.  For instance, all J\o{}rgensen subgroups of the Picard group $\mathrm{PSL}_2 (O_1)$ are found in \cite{Gonzalez-AcunaRamirez07}, and Theorem 3 of \cite{JorgensenKiikka75} solves Problem \ref{prb} in the Fuchsian case:

\begin{thm} \label{fuchjorg}
The only Fuchsian J\o{}rgensen groups are the $(2,3,q)$-triangle groups where $q \geq 7$ or $q = \infty$.
\end{thm}

This paper solves Problem \ref{prb} in two cases: torsion-free Kleinian groups and non-cocompact arithmetic Kleinian groups (recall that a finite-covolume Kleinian group is non-cocompact if and only if it contains parabolic elements; see, for instance, Section 1.2 of \cite{MaclachlanReid03}).

Identifying $\pslc$ with the orientation-preserving isometries of $\mathbb{H}^3$, we exploit in Section \ref{manifolds} the correspondence of torsion-free Kleinian groups $\Gamma$ with orientable hyperbolic 3-manifolds $\mathbb{H}^3/\Gamma$ (see, for instance, Section 1.3 of \cite{MaclachlanReid03}) to prove the following result.

\begin{thm} \label{fig8}
The figure-eight knot group is the only torsion-free J\o{}rgensen group.
\end{thm}

We then observe in Section \ref{los} that the main theorems of \cite{LiOichiSato04}, \cite{LiOichiSato05a}, and \cite{LiOichiSato05b} combine to find all J\o{}rgensen groups of parabolic type $(\theta, k)$, meaning that the group can be generated by the matrices
\[ A = \left(\begin{array}{cc} 1&1 \\ 0&1 \\ \end{array}\right) \text{ and }
B_{\theta, k} = \left(\begin{array}{rc} ke^{i\theta}&ik^2e^{i\theta}-ie^{-i\theta} \\ -ie^{i\theta}&ke^{i\theta} \\ \end{array}\right) \]
where $0 \leq \theta \leq 2\pi$ and $k \in \R$ (note our mild abuse of notation which will be continued tacitly throughout: we blur the distinction between elements of $\mathrm{PSL}_2\mathbb{C}$ and their lifts to $\mathrm{SL}_2\mathbb{C}$).  It is conjectured in \cite{LiOichiSato04}, \cite{LiOichiSato05a}, and \cite{LiOichiSato05b} that every J\o{}rgensen group of parabolic type is conjugate in $\pslc$ to a J\o{}rgensen group of parabolic type $(\theta, k)$.  If this conjecture were true, then Problem \ref{prb} would be solved for J\o{}rgensen groups of parabolic type.

In Section \ref{arith}, however, we disprove this conjecture by exhibiting four J\o{}rgensen groups of parabolic type that are not conjugate to any J\o{}rgensen group of parabolic type $(\theta, k)$ found in \cite{LiOichiSato04}, \cite{LiOichiSato05a}, or \cite{LiOichiSato05b}; these counterexamples are $\mathrm{PGL}_2 (O_3)$, $\mathrm{PSL}_2 (O_3)$, $\mathrm{PSL}_2 (O_7)$, and $\mathrm{PSL}_2 (O_{11})$, where $O_d$ is the ring of integers in $\mathbb{Q}(\sqrt{-d})$ for $d \in \mathbb{N}$. This is part of our second contribution to Problem \ref{prb}: we identify all non-cocompact arithmetic J\o{}rgensen groups (and hence all arithmetic J\o{}rgensen groups of parabolic type) by proving the following result.

\begin{thm} \label{parabcase}
There are exactly fourteen non-cocompact arithmetic J\o{}rgensen groups; they are: \\ $\mathrm{PGL}_2 (O_1)$, $\mathrm{PGL}_2 (O_3)$, $\mathrm{PSL}_2 (O_1)$, $\mathrm{PSL}_2 (O_2)$, $\mathrm{PSL}_2 (O_3)$, $\mathrm{PSL}_2 (O_7)$, $\mathrm{PSL}_2 (O_{11})$, two subgroups of index 6 and 8 respectively in $\mathrm{PGL}_2 (O_1)$, the unique subgroup of index 10 in $\mathrm{PSL}_2 (O_3)$, the figure-eight knot group, and three $\mathbb{Z}_2$-extensions of the figure-eight knot group.
\end{thm}

Since J\o{}rgensen groups are two-generator by definition, Theorem \ref{parabcase} can also be regarded as a solution to the following problem in the case of non-cocompact J\o{}rgensen groups.

\begin{prob} \label{prb2}
Identify all two-generator arithmetic Kleinian groups.
\end{prob}

Problem \ref{prb2} has also attracted much attention, including \cite{MaclachlanMartin99}, which proves that only finitely many arithmetic Kleinian groups can be generated by two elliptic elements; \cite{GehringMaclachlanMartin98}, which identifies all arithmetic Kleinian groups generated by two parabolic elements; and \cite{Conderetal02}, which identifies all two-generator arithmetic Kleinian groups with one generator parabolic and the other elliptic (see Theorem \ref{conder} for statement).

In proving Theorem \ref{parabcase}, we demonstrate that $\mathrm{PSL}_2 (O_1)$, $\mathrm{PSL}_2 (O_2)$, $\mathrm{PSL}_2 (O_7)$, $\mathrm{PSL}_2 (O_{11})$, the unique subgroup of index 10 in $\mathrm{PSL}_2 (O_3)$, a $\mathbb{Z}_2$-extension of the figure-eight knot group by an involution that conjugates each parabolic generator to its own inverse, and the two subgroups of index 6 and 8 respectively in $\mathrm{PGL}_2 (O_1)$ are two-generator arithmetic Kleinian groups with one generator parabolic and the other loxodromic. None of these were identified as two-generator arithmetic Kleinian groups in \cite{MaclachlanMartin99}, \cite{GehringMaclachlanMartin98}, or \cite{Conderetal02} and can therefore be regarded as a contribution to Problem \ref{prb2}.

We now turn our attention to J\o{}rgensen groups of elliptic type.  The following useful result is observed in \cite{GehringMartin89}.

\begin{thm} \label{subtriang}
Every J\o{}rgensen group of elliptic type contains a $(2,3,q)$-triangle group where $q \geq 7$.
\end{thm}

Exploiting this and the enumeration of all arithmetic triangle groups in \cite{Takeuchi77}, we establish in Section \ref{elliptic} a characterization of arithmetic J\o{}rgensen groups of elliptic type from which it follows that they are cocompact.

We conclude in Section \ref{bounds} by bounding $J(\Gamma)$ if $\Gamma$ is non-cocompact (e.g., a two-bridge knot or link group) and in Section \ref{comp} by computing $J(\Gamma)$ for several such groups, including that of the Whitehead link, whose J\o{}rgensen number was shown to be two in \cite{Sato04}; the others have not yet appeared in publication.

\section{Torsion-Free J\o{}rgensen Groups} \label{manifolds}

To prove Theorem \ref{fig8}, we first recall the notion of waist size for cusps in hyperbolic 3-manifolds introduced in \cite{Adams02}.  The \emph{waist size} $w(M, \mathcal{C})$ of a cusp $\mathcal{C}$ in an orientable hyperbolic 3-manifold $M = \mathbb{H}^3/\Gamma$ is the length of the shortest translation in $\mathcal{C}$ after expanding $\mathcal{C}$ until it first touches itself on its boundary, i.e., expanding $\mathcal{C}$ to a \emph{maximal cusp}.  Note that if $\Gamma$ contains parabolic elements, then $M$ is non-compact and contains one or more cusps. The following is a combination of Lemma 2.4 and Theorem 3.1, both from \cite{Adams02}.

\begin{thm} \label{adams}
The waist size of any cusp in an orientable hyperbolic 3-manifold is at least one, and the only hyperbolic 3-manifold with a cusp of waist size one is the figure-eight knot complement in $S^3$.
\end{thm}

We begin by bounding waist size 
(cf. Lemma 2.5 of \cite{BakerReid02}).

\begin{lem} \label{waistlem}
Let $M = \mathbb{H}^3/\Gamma$ be an orientable hyperbolic 3-manifold containing a cusp $\mathcal{C}$. Conjugate $\Gamma$ so that $\mathcal{C}$ lifts to a disjoint set of horoballs, one of which, $\mathcal{H}$, is based at $\infty$, and so that $\Gamma$ contains the parabolic element $A = \left(\begin{array}{cc} 1&1 \\ 0&1 \\ \end{array}\right)$.  Then \[ 1 \leq w(M, \mathcal{C}) \leq \inf \left\{|c| : \left(\begin{array}{cc} a&b \\ c&d \\ \end{array}\right) \in \Gamma \text{ and } c \neq 0 \right\}. \]
\end{lem}

\begin{proof}
After expanding $\mathcal{C}$ to a maximal cusp, its height $h$ is that of the expanded horoball $\mathcal{H}$. Let $T = \left(\begin{array}{cc} a&b \\ c&d \\ \end{array}\right) \in \Gamma$ with $c \neq 0$.  Then $T(\mathcal{H})$ is a horoball based at $\dfrac{a}{c}$ of diameter $\dfrac{1}{h|c|^2}$ whose interior is disjoint from that of $\mathcal{H}$.  Hence, $h \geq \dfrac{1}{h|c|^2}$, so $h \geq \dfrac{1}{|c|}$.  The translation length of $A$ at height $h$ is $\dfrac{1}{h} \leq |c|$, so $w(M,
\mathcal{C}) \leq |c|$.  Since $T$ was arbitrary,
\[ w(M, \mathcal{C}) \leq \inf \left\{|c| : \left(\begin{array}{cc} a&b \\ c&d \\ \end{array}\right) \in \Gamma \text{ and } c \neq 0 \right\}. \]
Theorem \ref{adams} establishes $1 \leq w(M, \mathcal{C})$.
\end{proof}

To prove Theorem \ref{fig8}, we define the \emph{generalized J\o{}rgensen number} of an arbitrary subgroup $\Gamma$ of $\pslc$ to be
\[ \widetilde{J}(\Gamma) = \inf\{J(X, Y) : \langle X, Y \rangle \leq \Gamma \text{ is discrete and non-elementary} \}. \]
Again, J\o{}rgensen's Inequality guarantees $\widetilde{J}(\Gamma) \geq 1$, and if $\Gamma$ is rank-two, non-elementary, and discrete, then $J(\Gamma) \geq \widetilde{J}(\Gamma)$ by definition.

\begin{thm} \label{part1}
The only orientable hyperbolic 3-manifold $\mathbb{H}^3/\Gamma$ with $\widetilde{J}(\Gamma) = 1$ is the figure-eight knot complement in $S^3$.
\end{thm}

\begin{proof}
As established in \cite{Riley75}, if $\mathbb{H}^3/\Gamma$ is the figure-eight knot complement in $S^3$, then $\Gamma$ is generated by the matrices
\[A = \left(\begin{array}{cc} 1&1 \\ 0&1 \\ \end{array}\right)
\text{ and }B = \left(\begin{array}{cc} 1&0 \\ e^{\pi i/3}&1 \\ \end{array}\right).\]
A simple computation shows that $J(A,B) = |e^{2\pi i/3}| = 1$, so $\widetilde{J}(\Gamma) = J(\Gamma) = 1$.

Conversely, let $M = \mathbb{H}^3/\Gamma$ be an orientable hyperbolic 3-manifold with $\widetilde{J}(\Gamma) = 1$. Then there exists a non-elementary subgroup $\langle A, B \rangle \leq \Gamma$ with $J(A, B) = 1$. Since $\Gamma$ is torsion-free, Theorem \ref{kiikka} implies $A$ is parabolic. Conjugate $\Gamma$ so that
\[ A = \left(\begin{array}{cc} 1&1 \\ 0&1 \\ \end{array}\right) \text{ and } B = \left(\begin{array}{cc} a&b \\ c&d \\ \end{array}\right). \]
Then $M$ contains a cusp $\mathcal{C}$ that lifts to a disjoint set of horoballs, one of which is based at $\infty$. By assumption, \[J(A,B) = |\tr^2 A - 4| + |\tr[A,B] - 2| = |c|^2 = 1,\] so $|c| = 1$, and $w(M,\mathcal{C}) = 1$ by Lemma \ref{waistlem}.  Hence, $\mathbb{H}^3/\Gamma$ is the figure-eight knot complement in $S^3$ by Theorem \ref{adams}.
\end{proof}

Since $J(\Gamma) \geq \widetilde{J}(\Gamma) \geq 1$, Theorem \ref{fig8} now follows easily:

\begin{cor}
The only orientable hyperbolic 3-manifold $\mathbb{H}^3/\Gamma$ with $J(\Gamma) = 1$ is the figure-eight knot complement in $S^3$.
\end{cor}

\section{J\o{}rgensen Groups of Parabolic Type $(\theta, k)$} \label{los}

Following \cite{LiOichiSato04}, \cite{LiOichiSato05a}, and \cite{LiOichiSato05b}, a group of parabolic type $(\theta, k)$ is denoted $G_{\theta, k}$ and is generated by the matrices
\[ A = \left(\begin{array}{cc} 1&1 \\ 0&1 \\ \end{array}\right) \text{ and }
B_{\theta, k} = \left(\begin{array}{rc} ke^{i\theta}&ik^2e^{i\theta}-ie^{-i\theta} \\ -ie^{i\theta}&ke^{i\theta} \\ \end{array}\right) \]
where $0 \leq \theta \leq 2\pi$ and $k \in \R$.  The main theorems of \cite{LiOichiSato04}, \cite{LiOichiSato05a}, and \cite{LiOichiSato05b} combine to find all J\o{}rgensen groups of parabolic type $(\theta, k)$, and the following conjecture is made.

\begin{con} \label{losconj}
Every J\o{}rgensen group of parabolic type is conjugate in $\pslc$ to a J\o{}rgensen group of parabolic type $(\theta, k)$.
\end{con}

Thus, if Conjecture \ref{losconj} were true, then Problem \ref{prb} is solved for J\o{}rgensen groups of parabolic type.  We will see in Section \ref{arith} that this is not the case by exhibiting four arithmetic J\o{}rgensen groups of parabolic type that are not conjugate to any J\o{}rgensen group of parabolic type $(\theta, k)$ found in \cite{LiOichiSato04}, \cite{LiOichiSato05a}, or \cite{LiOichiSato05b}; these counterexamples are $\mathrm{PGL}_2 (O_3)$, $\mathrm{PSL}_2 (O_3)$, $\mathrm{PSL}_2 (O_7)$, and $\mathrm{PSL}_2 (O_{11})$.  Since arithmetic Kleinian groups necessarily have finite covolume, we restrict our attention to finite-covolume J\o{}rgensen groups of parabolic type $(\theta, k)$ and state the main theorems of \cite{LiOichiSato04}, \cite{LiOichiSato05a}, and \cite{LiOichiSato05b}, together with Corollary 3.5, Lemma 3.6, and Lemma 3.8 of \cite{LiOichiSato05a} (cf. Corollary 6.3 and Lemma 6.4 of \cite{Sato00}), as follows.

\begin{thm} \label{losthm}
For $0 \leq \theta \leq \pi/2$ and $k \geq 0$, $G_{\theta, k}$ is a finite-covolume J\o{}rgensen group if and only if $(\theta, k)$ is one of the following pairs.\\
\begin{tabular}{llll}
$\bullet \; (\frac{\pi}{6},\frac{\sqrt{3}}{2}n)$ for $n \in \Z$&$\bullet \; (\frac{\pi}{4},\frac{1}{2})$&$\bullet \; (\frac{\pi}{4},1)$&$\bullet \; (\frac{\pi}{4},\frac{3}{2})$\\
$\bullet \; (\frac{\pi}{4},1+\frac{\sqrt{2}}{2})$&$\bullet \; (\frac{\pi}{4},\frac{5+\sqrt{5}}{4})$&$\bullet \; (\frac{\pi}{4},1+\frac{\sqrt{3}}{2})$&$\bullet \; (\frac{\pi}{3},\frac{\sqrt{3}}{2}n)$ for $n \in \Z$\\
$\bullet \; (\frac{\pi}{2},\frac{1}{2})$&$\bullet \; (\frac{\pi}{2},\frac{\sqrt{2}}{2})$&$\bullet \; (\frac{\pi}{2},\frac{1+\sqrt{5}}{4})$&$\bullet \; (\frac{\pi}{2},\frac{\sqrt{3}}{2})$\\
\end{tabular}

\noindent Furthermore,
\begin{itemize}
\item For $0 \leq \theta \leq \pi/2$ and $k \in \R$, $G_{\theta, k}$ is a J\o{}rgensen group if and only if $G_{\pi-\theta, k}$ is a J\o{}rgensen group.
\item For $0 \leq \theta \leq \pi$ and $k \in \R$, $G_{\pi+\theta, k} = G_{\theta, k}$.
\item For $0 \leq \theta \leq 2\pi$ and $k \in \R$, $G_{\theta, k}$ is a J\o{}rgensen group if and only if $G_{\theta, -k}$ is a J\o{}rgensen group.
\end{itemize}
\end{thm}

This completes the identification of all finite-covolume J\o{}rgensen groups of parabolic type $(\theta, k)$.  To determine which of these is arithmetic, we recall some basic facts about two-generator Kleinian groups and arithmeticity; these will also be useful in the next two sections.  The first is a combination of (3.25) and Lemmas 3.5.7, 3.5.8, and 8.5.2 in \cite{MaclachlanReid03}.

\begin{lem} \label{tracefields}
Let $\Gamma = \langle A,B \rangle$ be a non-elementary Kleinian group with $\tr A \neq 0$.
\begin{itemize}
\item The trace field is \[\Q(\tr\Gamma) = \Q(\tr A, \tr B, \tr AB).\]
\item If $\tr B = 0$, then the invariant trace field is \[k\Gamma = \mathbb{Q}(\tr^2 A, \tr [A,B]).\]
\item If $\tr B \neq 0$, then the invariant trace field is \[k\Gamma = \mathbb{Q}(\tr^2 A, \tr^2B, \tr A \tr B \tr AB).\]
\item If $\tr A$, $\tr B$, and $\tr AB$ are algebraic integers, then $\tr \Gamma$ consists of algebraic integers.
\end{itemize}
\end{lem}

The following characterization will enable us to determine when a finite-covolume non-cocompact Kleinian group is arithmetic (cf. Theorem 2.2 of \cite{GehringMaclachlanMartin98} and Theorem 8.2.3 of \cite{MaclachlanReid03}).

\begin{thm} \label{bianchi}
A finite-covolume non-cocompact Kleinian group $\Gamma$ is arithmetic if and only if $\tr \Gamma$ consists of algebraic integers and $k\Gamma = \Q(\sqrt{-d})$ for some $d \in \N$.
\end{thm}

This is equivalent to $\Gamma$ being commensurable with the Bianchi group $\mathrm{PSL}_2 (O_d)$, so it is useful to recall from Section 5 of \cite{Conderetal02} the presentations of several Bianchi groups (cf. \cite{GrunewaldSchwermer93} and \cite{Swan71}).

\begin{thm} \label{pres}
For $d \in \{1, 2, 3, 7, 11\}$, $\mathrm{PSL}_2 (O_d)$ is generated by the three matrices
\[ A = \left(\begin{array}{cc} 1&1 \\ 0&1 \\ \end{array}\right), \;  S = \left(\begin{array}{cr} 0&-1 \\ 1&0 \\ \end{array}\right), \text{ and }  T = \left(\begin{array}{cc} 1&\alpha \\ 0&1 \\ \end{array}\right)\] where $\alpha = \sqrt{-d}$ if $d \not \equiv 3 \mod 4$ and $\alpha = (1 + \sqrt{-d})/2$ if $d \equiv 3 \mod 4$.  The full presentations of these Bianchi groups are as follows.
\begin{itemize}
\item $\mathrm{PSL}_2 (O_1) = \langle A, S, T \; | \; S^2 = (AS)^3 = [A,T] = (T^2ST^{-1}S)^2 = $ \begin{flushright}$(TST^{-1}STS)^2 = (ATST^{-1}STS)^2 = 1\rangle$\end{flushright}
\item  $\mathrm{PSL}_2 (O_2) = \langle A, S, T \; | \; S^2 = (AS)^3 = [A,T] = [S,T]^2 = 1\rangle$
\item  $\mathrm{PSL}_2 (O_3) = \langle A, S, T \; | \; S^2 = (AS)^3 = [A,T] = (TST^{-1}AT^{-1}AS)^2 = $ \begin{flushright}$(TST^{-1}AS)^3 = A^{-1}T^{-1}SA^{-1}TSAT^{-1}SAT^{-1}SA^{-1}TS = 1\rangle$\end{flushright}
\item  $\mathrm{PSL}_2 (O_7) = \langle A, S, T \; | \; S^2 = (AS)^3 = [A,T] = (SAT^{-1}ST)^2 = 1\rangle$
\item  $\mathrm{PSL}_2 (O_{11}) = \langle A, S, T \; | \; S^2 = (AS)^3 = [A,T] = (SAT^{-1}ST)^3 = 1\rangle$
\end{itemize}
\end{thm}

We also recall Theorem 1.1 from \cite{Conderetal02}, which identifies all two-generator arithmetic Kleinian groups with one generator parabolic and the other elliptic.  We further note the invariant trace field of each.

\begin{thm} \label{conder}
Suppose $\Gamma = \langle A,B \rangle$ is an arithmetic Kleinian group with $A$ parabolic and $B$ elliptic.  Then $B$ has order 2, 3, 4, or 6, and there are exactly fourteen such groups:
\begin{itemize}
\item If $B$ has order 2, then there are six groups:
\begin{enumerate}
\item Two $\mathbb{Z}_2$-extensions of the figure eight knot group each with index 6 in $\mathrm{PSL}_2 (O_3)$; $k\Gamma = \Q(\sqrt{-3})$.
\item A $\mathbb{Z}_2$-extension of the Whitehead link group with $\Gamma \cap \mathrm{PSL}_2 (O_1)$ of index 2 in $\Gamma$ and 12 in $\mathrm{PSL}_2 (O_1)$; $k\Gamma = \Q(i)$.
\item Two $\mathbb{Z}_2$-extensions of the $6^2_2$ link group each with $\Gamma \cap \mathrm{PSL}_2 (O_3)$ of index 2 in $\Gamma$ and 24 in $\mathrm{PSL}_2 (O_3)$; $k\Gamma = \Q(\sqrt{-3})$.
\item A $\mathbb{Z}_2$-extension of the $6^2_3$ link group with $\Gamma \cap \mathrm{PSL}_2 (O_7)$ of index 2 in $\Gamma$ and 12 in $\mathrm{PSL}_2 (O_7)$; $k\Gamma = \Q(\sqrt{-7})$.
\end{enumerate}
\item If $B$ has order 3, then there are three groups:
\begin{enumerate}
\item $[\mathrm{PSL}_2 (O_1) : \Gamma] = 8$; $k\Gamma = \Q(i)$.
\item $\Gamma = \mathrm{PSL}_2 (O_3)$; $k\Gamma = \Q(\sqrt{-3})$.
\item $[\mathrm{PSL}_2 (O_7) : \Gamma] = 2$; $k\Gamma = \Q(\sqrt{-7})$.
\end{enumerate}
\item If $B$ has order 4, then there are three groups:
\begin{enumerate}
\item $\Gamma = \mathrm{PGL}_2 (O_1)$; $k\Gamma = \Q(i)$.
\item $\Gamma \cap \mathrm{PSL}_2 (O_2)$ has index 4 in $\Gamma$ and 24 in $\mathrm{PSL}_2 (O_2)$; $k\Gamma = \Q(\sqrt{-2})$.
\item $\Gamma \cap \mathrm{PSL}_2 (O_3)$ has index 2 in $\Gamma$ and 30 in $\mathrm{PSL}_2 (O_3)$; $k\Gamma = \Q(\sqrt{-3})$.
\end{enumerate}
\item If $B$ has order 6, then there are two groups:
\begin{enumerate}
\item $\Gamma = \mathrm{PGL}_2 (O_3)$; $k\Gamma = \Q(\sqrt{-3})$.
\item $\Gamma \cap \mathrm{PSL}_2 (O_{15})$ has index 6 in $\Gamma$ and 6 in $\mathrm{PSL}_2 (O_{15})$; $k\Gamma = \Q(\sqrt{-15})$.
\end{enumerate}
\end{itemize}
\end{thm}

We are now prepared to identify which J\o{}rgensen groups of parabolic type $(\theta, k)$ in Theorem \ref{losthm} are arithmetic.  We conjugate $G_{\theta, k}$ by $\left(\begin{array}{cc} 1&-ik \\ 0&1 \\ \end{array}\right)$ and henceforth regard $G_{\theta, k}$ as the group generated by the matrices \[A = \left(\begin{array}{cc} 1&1 \\ 0&1 \\ \end{array}\right) \text{ and } B_{\theta, k} = \left(\begin{array}{cc} 0&-ie^{-i\theta} \\ -ie^{i\theta}&2ke^{i\theta} \\ \end{array}\right)\] where $0 \leq \theta \leq 2\pi$ and $k \in \R$.

\begin{pro} \label{losid}
For $0 \leq \theta \leq \pi/2$ and $k \geq 0$, $G_{\theta, k}$ is an arithmetic J\o{}rgensen group if and only if $(\theta, k)$ is one of the following pairs.
\begin{itemize}
\item $(\frac{\pi}{6},\frac{\sqrt{3}}{2}n)$ for $n \in \Z$, in which case $G_{\theta, k}$ is the figure-eight knot group if $n$ is odd and a $\Z_2$-extension of the figure-eight group by an involution that conjugates one parabolic generator to the other if $n$ is even.
\item $(\frac{\pi}{4},\frac{1}{2})$, in which case $G_{\theta, k} = \mathrm{PGL}_2 (O_1)$.
\item $(\frac{\pi}{4},1)$, in which case $G_{\theta, k}$ is a subgroup of index 8 in $\mathrm{PGL}_2 (O_1)$.
\item $(\frac{\pi}{4},\frac{3}{2})$, in which case $G_{\theta, k}$ is a subgroup of index 6 in $\mathrm{PGL}_2 (O_1)$.
\item $(\frac{\pi}{3},\frac{\sqrt{3}}{2}n)$ for $n \in \Z$, in which case $G_{\theta, k}$ is a $\Z_2$-extension of the figure-eight group by an involution that conjugates one parabolic generator to the other if $n$ is even and by an involution that conjugates each parabolic generator to its own inverse if $n$ is odd.
\item $(\frac{\pi}{2},\frac{1}{2})$, in which case $G_{\theta, k} = \mathrm{PSL}_2 (O_1)$.
\item $(\frac{\pi}{2},\frac{\sqrt{2}}{2})$, in which case $G_{\theta, k} = \mathrm{PSL}_2 (O_2)$.
\item $(\frac{\pi}{2},\frac{\sqrt{3}}{2})$, in which case $G_{\theta, k}$ is the unique subgroup of index 10 in $\mathrm{PSL}_2 (O_3)$.
\end{itemize}
\end{pro}

\begin{proof}
Let $G_{\theta, k}$ be a finite-covolume J\o{}rgensen group with $0 \leq \theta \leq \pi/2$ and $k \geq 0$.  Then $(\theta, k)$ must be one of the pairs listed in Theorem \ref{losthm}.  By Lemma \ref{tracefields}, $\tr G_{\theta, k}$ consists of algebraic integers for each of these pairs, but if $(\theta, k) = (\frac{\pi}{4}, 1+\frac{\sqrt{2}}{2})$, $(\frac{\pi}{4}, \frac{5+\sqrt{5}}{4})$, $(\frac{\pi}{4}, 1+\frac{\sqrt{3}}{2})$, or $(\frac{\pi}{2}, \frac{5+\sqrt{5}}{4})$, then $kG_{\theta, k} = \Q(ki, e^{2i\theta})$ is not a quadratic imaginary extension of $\Q$ and hence $G_{\theta, k}$ is not arithmetic by Theorem \ref{bianchi}.

For each of the remaining cases, we first use Lemma \ref{tracefields} to see that $kG_{\theta, k}$ is a quadratic imaginary extension of $\Q$ (and hence $G_{\theta, k}$ is arithmetic by Theorem \ref{bianchi}), and then we identify $G_{\theta, k}$ using the theorems above.

\textbf{Case 1.} Suppose $(\theta, k) = (\frac{\pi}{6}, \frac{\sqrt{3}}{2}n)$ for some integer $n$.  Then $kG_{\theta, k} = \Q(\sqrt{-3})$. Let
\begin{eqnarray*}
C &=& B_{\theta,k}AB_{\theta,k}^{-1} = \left(\begin{array}{cc} 1&0 \\ e^{\pi i/3}&1 \\ \end{array}\right) \text{ and }\\
T &=& C^{-1}ACA^{-2}CAC^{-1} = \left(\begin{array}{cc} 1&2\sqrt{3}i \\ 0&1 \\ \end{array}\right).
\end{eqnarray*}
We have already seen that $\langle A,C \rangle$ is the figure-eight knot group.

If $n$ is even, then let
\[D = B_{\theta,k}T^{-n/2} = \left(\begin{array}{cc} 0&-e^{\pi i/3} \\ e^{-\pi i/3}&0 \\ \end{array}\right).\]
Then $B_{\theta,k} = DT^{n/2}$, $T \in \langle A,C \rangle$, and $DAD^{-1} = C$, so $G_{\theta, k} = \langle A, D \rangle$, which is a $\Z_2$-extension of the figure eight knot group by an involution that conjugates one generator to the other.  Note also that $\Q(\tr G_{\theta, k}) = \Q(e^{\pi i/3})$ by Lemma \ref{tracefields}.

If $n$ is odd, then observe that
\[B_{\theta,k}^{-1}A^{-1}CAC^{-1}T^{(n-1)/2} = \left(\begin{array}{cc} 1&0 \\ 0&1 \\ \end{array}\right),\]
so $B_{\theta,k} \in \langle A, C \rangle$ since $T \in \langle A, C \rangle$.  Thus $G_{\theta, k} = \langle A, C \rangle$, the figure-eight knot group.

\textbf{Case 2.} Suppose $(\theta, k) = (\frac{\pi}{4}, \frac{1}{2})$.  Then $kG_{\theta, k} = \Q(i)$ and $G_{\theta, k} = \langle A, B_{\theta,k}A \rangle$.  Since $B_{\theta,k}A$ has order 4, we conclude that $G_{\theta, k} = \mathrm{PGL}_2 (O_1)$ by Theorem \ref{conder}.

\textbf{Case 3.} Suppose $(\theta, k) = (\frac{\pi}{4}, 1)$.  Then $kG_{\theta, k} = \Q(i)$. By Section 6.4.2 of \cite{LiOichiSato05a}, \[G_{\frac{\pi}{4}, \frac{1}{2}} = \langle S,T,U \mid U^2 = S^4 = T^4 = (US)^2 = (U^{-1}T)^3 = (TS)^2 =1 \rangle\]
where
\begin{eqnarray*}
S &=& AB_{\frac{\pi}{4}, \frac{1}{2}}^2AB_{\frac{\pi}{4}, \frac{1}{2}}^{-1},\\
T &=& A^2B_{\frac{\pi}{4}, \frac{1}{2}}^2AB_{\frac{\pi}{4}, \frac{1}{2}}^{-1},\text{ and }\\
U &=& AB_{\frac{\pi}{4}, \frac{1}{2}}^2A^{-1}B_{\frac{\pi}{4}, \frac{1}{2}}^{-2}A^{-1}B_{\frac{\pi}{4}, \frac{1}{2}}.
\end{eqnarray*}
Then $A = TS^{-1}$ and $B_{\frac{\pi}{4}, 1} = SUTS^2T^{-1}S^{-1}$, so, using Magma (\cite{Magma}), we find that $G_{\frac{\pi}{4}, 1}$ is a subgroup of index 8 in $G_{\frac{\pi}{4}, \frac{1}{2}} = \mathrm{PGL}_2 (O_1)$.

\textbf{Case 4.} Suppose $(\theta, k) = (\frac{\pi}{4}, \frac{3}{2})$.  Then $kG_{\theta, k} = \Q(i)$.  Using the presentation of $G_{\frac{\pi}{4}, \frac{1}{2}}$ above, we see that $A = TS^{-1}$ and \[B_{\frac{\pi}{4}, \frac{3}{2}} = SUTS^2T^{-1}ST^{-1}S^{-1}.\] Using Magma (\cite{Magma}), we find that $G_{\frac{\pi}{4}, \frac{3}{2}}$ is a subgroup of index 6 in $G_{\frac{\pi}{4}, \frac{1}{2}} = \mathrm{PGL}_2 (O_1)$.

\textbf{Case 5.} Suppose $(\theta, k) = (\frac{\pi}{3}, \frac{\sqrt{3}}{2}n)$ for some integer $n$. Then $kG_{\theta, k} = \Q(\sqrt{-3})$. Let
\begin{eqnarray*}
C &=& B_{\theta,k}A^{-1}B_{\theta,k}^{-1} = \left(\begin{array}{cc} 1&0 \\ e^{-\pi i/3}&1 \\ \end{array}\right) \text{ and }\\
T &=& CA^{-1}C^{-1}A^2C^{-1}A^{-1}C = \left(\begin{array}{cc} 1&2\sqrt{3}i \\ 0&1 \\ \end{array}\right).
\end{eqnarray*}
Note that $\langle A,C \rangle$ is again the figure-eight knot group.

If $n$ is even, then let
\[D = B_{\theta,k}T^{-n/2} = \left(\begin{array}{cc} 0&-e^{\pi i/6} \\ e^{-\pi i/6}&0 \\ \end{array}\right).\]
Then $B_{\theta,k} = DT^{n/2}$, $T \in \langle A,C \rangle$, and $DAD^{-1} = C$, so $G_{\theta, k} = \langle A, D \rangle$, which is a $\Z_2$-extension of the figure eight knot group by an involution that conjugates one generator to the other, but here $\Q(\tr G_{\theta, k}) = \Q(e^{\pi i/6})$ by Lemma \ref{tracefields}, so this $\Z_2$-extension of the figure eight knot group is not conjugate to the one found in Case 1 since the trace field is a conjugacy invariant of finite-covolume Kleinian groups (see, for instance, Section 3.1 of \cite{MaclachlanReid03}).

If $n$ is odd, then let
\[D = B_{\theta,k}^{-1}A^{-1}CAC^{-1}T^{-(n-1)/2} = \left(\begin{array}{cr} i&0 \\ 0&-i \\ \end{array}\right),\]
so $B_{\theta,k} \in \langle A, C, D \rangle$ since $T \in \langle A, C \rangle$.  As $D$ conjugates $A$ and $C$ to their own inverses, we conclude that $G_{\theta, k} = \langle A, C, D \rangle$ is a $\Z_2$-extension of the figure eight knot group by an involution that conjugates each parabolic generator to its own inverse.

\textbf{Case 6.} Suppose $(\theta, k) = (\frac{\pi}{2}, \frac{1}{2})$.  Then $kG_{\theta, k} = \Q(i)$. Let
\begin{eqnarray*}
S &=& A^{-1}B_{\theta,k}A^{-1}B_{\theta,k}^{-1}A^{-1} = \left(\begin{array}{cr} 0&-1 \\ 1&0 \\ \end{array}\right) \text{ and }\\
T &=& AB_{\theta,k}AB_{\theta,k}^{-1}AB_{\theta,k} = \left(\begin{array}{cc} 1&i \\ 0&1 \\ \end{array}\right).
\end{eqnarray*}
Then $ST = B_{\theta,k}$, so $G_{\theta, k} = \langle A, S, T \rangle = \mathrm{PSL}_2 (O_1)$ by Theorem \ref{pres}.

\textbf{Case 7.} Suppose $(\theta, k) = (\frac{\pi}{2}, \frac{\sqrt{2}}{2})$.  Then $kG_{\theta, k} = \Q(\sqrt{-2})$. Let
\begin{eqnarray*}
S &=& A^{-1}B_{\theta,k}A^{-1}B_{\theta,k}^{-1}A^{-1} = \left(\begin{array}{cr} 0&-1 \\ 1&0 \\ \end{array}\right) \text{ and }\\
T &=& AB_{\theta,k}AB_{\theta,k}^{-1}AB_{\theta,k} = \left(\begin{array}{cc} 1&\sqrt{-2} \\ 0&1 \\ \end{array}\right).
\end{eqnarray*}
Then $ST = B_{\theta,k}$, so $G_{\theta, k} = \langle A, S, T \rangle = \mathrm{PSL}_2 (O_2)$ by Theorem \ref{pres}.

\textbf{Case 8.} Suppose $(\theta, k) = (\frac{\pi}{2}, \frac{\sqrt{3}}{2})$.  Then $kG_{\theta, k} = \Q(\sqrt{-3})$. Let
\begin{eqnarray*}
S &=& A^{-1}B_{\theta,k}A^{-1}B_{\theta,k}^{-1}A^{-1} = \left(\begin{array}{cr} 0&-1 \\ 1&0 \\ \end{array}\right) \text{ and }\\
T &=& AB_{\theta,k}AB_{\theta,k}^{-1}AB_{\theta,k} = \left(\begin{array}{cc} 1&\sqrt{-3} \\ 0&1 \\ \end{array}\right).
\end{eqnarray*}
Then $ST = B_{\theta,k}$, so $G_{\theta, k} = \langle A, S, T \rangle$.  Using Theorem \ref{pres} and Magma (\cite{Magma}), we conclude that $\Gamma$ is the unique subgroup of index 10 in the Bianchi group $\mathrm{PSL}_2 (O_3)$ (cf. Section 6 of \cite{GrunewaldSchwermer93}).
\end{proof}

\section{Arithmetic J\o{}rgensen Groups of Parabolic Type} \label{arith}

To prove Theorem \ref{parabcase}, we first identify all arithmetic J\o{}rgensen groups of parabolic type.  We begin by stating Poincar\'{e}'s Polyhedron Theorem (see Section IV.H of \cite{Maskit88} for notation and terminology), which we will use to construct fundamental polyhedra for arithmetic J\o{}rgensen groups of parabolic type not treated in \cite{LiOichiSato04}, \cite{LiOichiSato05a}, or \cite{LiOichiSato05b}.

\newtheorem*{ppt}{Poincar\'{e}'s Polyhedron Theorem}
\begin{ppt}
Let $P$ be a polyhedron with side pairing transformations satisfying the following conditions (1) through (5).  Then the group $\Gamma$ generated by the side pairing transformations is discrete, $P$ is a fundamental polyhedron for $\Gamma$, and the reflection relations and cycle relations form a complete set of relations for $\Gamma$.
\begin{enumerate}
\item For each side $s$ of $P$, there is a side $s'$ and an element $g_s \in \Gamma$ such that $g_s(s) = s'$, $g_{s'} = g_s^{-1}$, and $g_s(P) \cap P = \emptyset$.
\item For every point $z \in P^*$, $p^{-1}(z)$ is a finite set, where $P^*$ is the space of equivalence classes so that the projection $p: \overline{P} \rightarrow P^*$ is continuous and open.
\item For each edge $e$, there is a positive integer $t$ so that $h^t = 1$, where $h$ is the cycle transformation at $e$.
\item If $\{e_1, e_2, \dots, e_m\}$ is any cycle of edges of $P$, $\alpha(e_k)$ the angle measured from inside $P$ at the edge $e_k$, and $q$ the smallest positive integer such that $h^q=1$, where $h$ is the cycle transformation at $e_k$, then \[ \sum_{k=1}^m \alpha(e_k) = 2\pi/q.\]
\item $P^*$ is complete.
\end{enumerate}
\end{ppt}

We also recall the notation of \cite{LiOichiSato04}, \cite{LiOichiSato05a}, and \cite{LiOichiSato05b}, which we will borrow when using Poincar\'{e}'s Polyhedron Theorem.  Let $F_X$ and $F_X^{-1}$ be two sides of a polyhedron $P$ such that $F_X$ is mapped onto $F_X^{-1}$ by the side pairing transformation $X$.  Denote by $e_{(m,n),\theta}$ the $n$th edge of the $m$th cycle transformation such that the angle measured from the polyhedron $P$ at the edge is $\theta$.  The following diagram represents the $m$th cycle transformation
\[ e_{(m,1),\theta_1} \stackrel{X_1}{\longrightarrow} e_{(m,2),\theta_2} \stackrel{X_2}{\longrightarrow} \cdots \stackrel{X_{m-1}}{\longrightarrow} e_{(m,n),\theta_n} \stackrel{X_n}{\longrightarrow} \circlearrowleft^p_\theta \]
where the initial edge $e_{(m,1),\theta_1}$ is mapped to the second edge $e_{(m,2),\theta_2}$ by the side pairing transformation $X_1$, then the edge $e_{(m,2),\theta_2}$ is mapped to the edge $e_{(m,3),\theta_3}$ by the side pairing transformation $X_2$, and so on.  
The final edge $e_{(m,n),\theta_n}$ is mapped to the initial edge $e_{(m,1),\theta_1}$ by the side pairing transformation $X_n$ and the sum of all angles at the edges in this sequence is equal to $\theta$.  The cycle transformation $X_n X_{n-1} \cdots X_1$ is either the identity transformation, in which case $p=1$, or an elliptic transformation of order $p$.

We are now prepared to identify all arithmetic J\o{}rgensen groups of parabolic type by proving the following theorem.  Throughout the proof, we note whenever a group of parabolic type $(\theta, k)$ is encountered and refer to its identification in Proposition  \ref{losid}.  To complete the proof of Theorem \ref{parabcase}, we will show in Section \ref{elliptic} that arithmetic J\o{}rgensen groups of elliptic type are cocompact and hence the fourteen arithmetic J\o{}rgensen groups of parabolic type listed below are precisely all non-cocompact arithmetic J\o{}rgensen groups.

\begin{thm} \label{parabtype}
There are exactly fourteen arithmetic J\o{}rgensen groups of parabolic type; they are: \\ $\mathrm{PGL}_2 (O_1)$, $\mathrm{PGL}_2 (O_3)$, $\mathrm{PSL}_2 (O_1)$, $\mathrm{PSL}_2 (O_2)$, $\mathrm{PSL}_2 (O_3)$, $\mathrm{PSL}_2 (O_7)$, $\mathrm{PSL}_2 (O_{11})$, two subgroups of index 6 and 8 respectively in $\mathrm{PGL}_2 (O_1)$, the unique subgroup of index 10 in $\mathrm{PSL}_2 (O_3)$, the figure-eight knot group, and three $\mathbb{Z}_2$-extensions of the figure-eight knot group.
\end{thm}

\begin{proof}
Let $\Gamma$ be an arithmetic J\o{}rgensen group of parabolic type.  Then $\Gamma = \langle A,B \rangle$ with $J(A,B)=1$ and $A$ parabolic.  Conjugate $\Gamma$ so that
\[A = \left(\begin{array}{cc} 1&1 \\ 0&1 \\ \end{array}\right) \text{ and } B = \left(\begin{array}{cc} 0&-1/\sigma \\ \sigma&\lambda \\ \end{array}\right) \text{ where }\sigma, \lambda \in \C, \sigma \neq 0\]
(cf. proof of Theorem 3.3 in \cite{Conderetal02}).  Working in $\pslc$, we may replace $B$ with $-B$ (and $\lambda$ with $-\lambda$).  Hence we assume $\re \, \sigma > 0$ or $\re \, \sigma = 0$ and $\im \, \sigma > 0$.

By Lemma \ref{tracefields} and Theorem \ref{bianchi}, $k\Gamma = \Q(\sqrt{-d})$ and $\tr[A,B]-2 = \sigma^2 \in O_d$ for some $d \in \N$.  Since $J(A,B) = |\sigma^2| = 1$, we see that $\sigma^2$ is a unit in $O_d$, so, given our assumption on $\sigma$, we have $\sigma \in \{1, i, e^{\pm \pi i/4}\}$ if $d=1$, $\sigma \in \{1, i, e^{\pm \pi i/6}, e^{\pm \pi i/3}\}$ if $d=3$, and $\sigma \in \{1, i\}$ otherwise.  We further note that $\tr B \tr AB - \tr^2 B \in O_d$ by Lemma \ref{tracefields} and Theorem \ref{bianchi}.

\textbf{Case 1.} Suppose $\sigma = 1$.  If $\tr B = \lambda = 0$, then $\Gamma = \mathrm{PSL}_2 \Z$, which is not an arithmetic Kleinian group, so $\tr B = \lambda \neq 0$ and $\tr B \tr AB - \tr^2 B = \lambda \in O_d$.  Hence, if $d \not \equiv 3 \mod 4$, then $\lambda = m + n\sqrt{-d}$ for some integers $m$ and $n$.  Replacing $B$ with $BA^{-m}$ (since $\langle A, BA^{-m} \rangle = \langle A, B \rangle$), assume $\lambda = n\sqrt{-d}$.  Then $\Gamma = G_{\theta, k}$ with $(\theta, k) = (\frac{\pi}{2}, \frac{\sqrt{d}}{2}n)$, so, by Proposition  \ref{losid}, $(\theta, k) = (\frac{\pi}{2}, \frac{1}{2})$ or $(\frac{\pi}{2}, \frac{\sqrt{2}}{2})$, which corresponds to $\mathrm{PSL}_2 (O_1)$ or $\mathrm{PSL}_2 (O_2)$ respectively.

Now suppose $d \equiv 3 \mod 4$. Then $\lambda = m + n\frac{1+\sqrt{-d}}{2}$ for some integers $m$ and $n$.  If $n$ is even, replace $B$ with $BA^{-m-n/2}$ and assume $\lambda = \frac{n\sqrt{-d}}{2}$.  Then $\Gamma = G_{\theta, k}$ with $(\theta, k) = (\frac{\pi}{2}, \frac{\sqrt{d}}{4}n)$, so, by Proposition  \ref{losid}, $(\theta, k) = (\frac{\pi}{2}, \frac{\sqrt{3}}{2})$, which corresponds to the unique subgroup of index 10 in $\mathrm{PSL}_2 (O_3)$.

If $n$ is odd, replace $B$ with $BA^{-m-\frac{n-1}{2}}$ and assume $\lambda = \frac{1+n\sqrt{-d}}{2}$.  Let
\begin{eqnarray*}
S &=& A^{-1}BA^{-1}B^{-1}A^{-1} = \left(\begin{array}{cr} 0&-1 \\ 1&0 \\ \end{array}\right),\\
T &=& S^{-1}B = \left(\begin{array}{cc} 1&\frac{1+n\sqrt{-d}}{2} \\ 0&1 \\ \end{array}\right), \text{ and }\\
U &=& TA^{-1} = \left(\begin{array}{cc} 1&\frac{-1+n\sqrt{-d}}{2} \\ 0&1 \\ \end{array}\right).
\end{eqnarray*}
Then $ST = B$ and $UA = T$, so $\Gamma = \langle A, S, T \rangle = \langle A, S, U \rangle$.  If $n<0$, replace $T$ with $U^{-1}$, so assume $n>0$.

If $n=1$ and $d = 3$, 7, or 11, then Theorem \ref{pres} identifies $\Gamma$ as $\mathrm{PSL}_2 (O_3)$, $\mathrm{PSL}_2 (O_7)$, or $\mathrm{PSL}_2 (O_{11})$ respectively. Thus, these Bianchi groups are arithmetic J\o{}rgensen groups of parabolic type but are not identified as such in Proposition  \ref{losid}.  Hence, $\mathrm{PSL}_2 (O_3)$, $\mathrm{PSL}_2 (O_7)$, and $\mathrm{PSL}_2 (O_{11})$ are not of the form $G_{\theta, k}$.

Now assume $n>1$ or $d>11$ and consider the polygon $P$ in Figure \ref{subcase122} with sides
\[F_A = \left\{(x,y,t) \in \mathbb{H}^3 \mid x = -\frac{1}{2}, \; \frac{1-n^2d}{4n\sqrt{d}} \leq y \leq \frac{n^2d-1}{4n\sqrt{d}}, \; x^2 + y^2 + t^2 \geq 1 \right\}\]
\[F_{A^{-1}} = \left\{(x,y,t) \in \mathbb{H}^3 \mid x = \frac{1}{2}, \; \frac{1-n^2d}{4n\sqrt{d}} \leq y \leq \frac{n^2d-1}{4n\sqrt{d}}, \; x^2 + y^2 + t^2 \geq 1 \right\}\]
\[F_S = \left\{(x,y,t) \in \mathbb{H}^3 \mid -\frac{1}{2} \leq x \leq 0, \; x^2 + y^2 + t^2 = 1 \right\}\]
\[F_{S^{-1}} = \left\{(x,y,t) \in \mathbb{H}^3 \mid 0 \leq x \leq \frac{1}{2}, \; x^2 + y^2 + t^2 = 1 \right\}\]
\[F_T = \left\{(x,y,t) \in \mathbb{H}^3 \mid -\frac{1}{2} \leq x \leq 0, \; y = \frac{1-n^2d-4x}{4n\sqrt{d}} \right\}\]
\[F_{T^{-1}} = \left\{(x,y,t) \in \mathbb{H}^3 \mid 0 \leq x \leq \frac{1}{2}, \; y = \frac{n^2d-1-4x}{4n\sqrt{d}} \right\}\]
\[F_U = \left\{(x,y,t) \in \mathbb{H}^3 \mid 0 \leq x \leq \frac{1}{2}, \; y = \frac{1+n^2d+4x}{4n\sqrt{d}} \right\}\]
\[F_{U^{-1}} = \left\{(x,y,t) \in \mathbb{H}^3 \mid -\frac{1}{2} \leq x \leq 0, \; y = \frac{n^2d+1+4x}{4n\sqrt{d}} \right\}\]
and edges

\begin{tabular}{lll}
$\bullet \; e_{(1),\pi} = F_S \cap F_{S^{-1}}$&$\bullet \; e_{(3,1),\theta} = F_A \cap F_{U^{-1}}$&$\bullet \; e_{(4,1),\theta} = F_A \cap F_T$\\
$\bullet \; e_{(2,1),\pi/3} = F_A \cap F_S$&$\bullet \; e_{(3,2),\theta} = F_{A^{-1}} \cap F_{T^{-1}}$&$\bullet \; e_{(4,2),\theta} = F_{A^{-1}} \cap F_U$\\
$\bullet \; e_{(2,2),\pi/3} = F_{A^{-1}} \cap F_{S^{-1}}$&$\bullet \; e_{(3,3),\phi} = F_T \cap F_U$&$\bullet \; e_{(4,3),\phi} = F_{T^{-1}} \cap F_{U^{-1}}$\\
\end{tabular}

\begin{figure}[htp]
\begin{center}
\input{subcase122.pstex_t}
\end{center}
\caption[A fundamental polyhedron in Subcase 1.2.2 of Theorem 3.16]{The polygon $P$ viewed from the point at $\infty$ above.  Solid lines denote edges, solid dots denote edges parallel to the $t$-axis, dashed lines are on the boundary of $\mathbb{H}^3$, and dotted lines denote the $x$- and $y$-axes in the boundary of $\mathbb{H}^3$.  The sides $F_S$ and $F_{S^{-1}}$ are on the unit sphere centered at the origin; all other sides are perpendicular to the boundary of $\mathbb{H}^3$.}
\label{subcase122}
\end{figure}

Then $X(F_X) = F_{X^{-1}}$ for $X \in \{A,S,T,U\}$ are side pairing transformations with the following four cycle transformations:
\begin{enumerate}
\item $e_{(1),\pi} \stackrel{S}{\longrightarrow} \circlearrowleft^2_\pi$
\item $e_{(2,1),\pi/3} \stackrel{A}{\longrightarrow} e_{(2,2),\pi/3} \stackrel{S^{-1}}{\longrightarrow}  \circlearrowleft^3_{2\pi/3}$
\item $e_{(3,1),\theta} \stackrel{A}{\longrightarrow} e_{(3,2),\theta} \stackrel{T^{-1}}{\longrightarrow} e_{(3,3),\phi} \stackrel{U}{\longrightarrow} \circlearrowleft^1_{2\pi}$
\item $e_{(4,1),\theta} \stackrel{A}{\longrightarrow} e_{(4,2),\theta} \stackrel{U}{\longrightarrow} e_{(3,3),\phi} \stackrel{T^{-1}}{\longrightarrow} \circlearrowleft^1_{2\pi}$
\end{enumerate}
Therefore, by Poincar\'{e}'s Polyhedron Theorem, $\Gamma$ has the presentation
\[ \Gamma = \langle A,S,T,U \mid S^2 = (S^{-1}A)^3 = UT^{-1}A = T^{-1}UA = 1 \rangle, \]
and $P$ is the fundamental polyhedron for $\Gamma$. Clearly $P$ has infinite volume, so $\Gamma$ is not arithmetic if $n>1$ or $d>11$.

\textbf{Case 2.} Suppose $\sigma = i$. If $\tr B = \lambda = 0$, then Lemma \ref{tracefields} yields $k\Gamma = \Q$, a contradiction, so $\tr B = \lambda \neq 0$ and $\tr B \tr AB - \tr^2 B = \lambda i \in O_d$. Hence, if $d \not \equiv 3 \mod 4$, then $\lambda i = m + n\sqrt{d}i$ for some integers $m$ and $n$, so $\lambda = n\sqrt{d}-mi$.  Replacing $B$ with $BA^m$, assume $\lambda = n\sqrt{d}$.  Then $\Gamma = G_{\pi+\theta, k} = G_{\theta, k}$ with $(\theta, k) = (0, \frac{n\sqrt{d}}{2})$, so $\Gamma$ is not a J\o{}rgensen group by Proposition  \ref{losid}.

Now suppose $d \equiv 3 \mod 4$. Then $\lambda i = m + n\frac{1+\sqrt{d}i}{2}$ for some integers $m$ and $n$, so $\lambda = \frac{n\sqrt{d}}{2}-(m+n/2)i$.  If $n$ is even, replace $B$ with $BA^{m+n/2}$ and assume $\lambda = \frac{n\sqrt{d}}{2}$.  Then $\Gamma = G_{\pi+\theta, k} = G_{\theta, k}$ with $(\theta, k) = (0, \frac{n\sqrt{d}}{4})$, so $\Gamma$ is not a J\o{}rgensen group by Proposition  \ref{losid}.

If $n$ is odd, replace $B$ with $BA^{m + \frac{n+1}{2}}$ and assume $\lambda = \frac{n\sqrt{d} + i}{2}$.  Let
\begin{eqnarray*}
S &=& A^{-1}BAB^{-1}A^{-1} = \left(\begin{array}{cr} 0&-1 \\ 1&0 \\ \end{array}\right),\\
T &=& S^{-1}B =  \left(\begin{array}{cc} i&\frac{n\sqrt{d} + i}{2} \\ 0&-i \\ \end{array}\right), \text{ and}\\
U &=& TA^{-1} = \left(\begin{array}{cc} i&\frac{n\sqrt{d} - i}{2} \\ 0&-i \\ \end{array}\right).
\end{eqnarray*}
Then $ST = B$ and $UA = T$, so $\Gamma = \langle A, S, T \rangle = \langle A, S, U \rangle$.  If $n<0$, replace $T$ with $U^{-1}$, so assume $n>0$.

If $n>1$ or $d>11$, then consider the polygon $P$ in Figure \ref{subcase222a} with sides
\[F_A = \left\{(x,y,t) \in \mathbb{H}^3 \mid x = -\frac{1}{2}, \; y \leq \frac{n^2d-1}{4n\sqrt{d}}, \; x^2 + y^2 + t^2 \geq 1 \right\}\]
\[F_{A^{-1}} = \left\{(x,y,t) \in \mathbb{H}^3 \mid x = \frac{1}{2}, \; y \leq \frac{n^2d-1}{4n\sqrt{d}}, \; x^2 + y^2 + t^2 \geq 1 \right\}\]
\[F_S = \left\{(x,y,t) \in \mathbb{H}^3 \mid -\frac{1}{2} \leq x \leq 0, \; x^2 + y^2 + t^2 = 1 \right\}\]
\[F_{S^{-1}} = \left\{(x,y,t) \in \mathbb{H}^3 \mid 0 \leq x \leq \frac{1}{2}, \; x^2 + y^2 + t^2 = 1 \right\}\]
\[F_T = \left\{(x,y,t) \in \mathbb{H}^3 \mid -\frac{1}{2} \leq x \leq -\frac{1}{4}, \; y = \frac{1+n^2d+4x}{4n\sqrt{d}} \right\}\]
\[F_{T^{-1}} = \left\{(x,y,t) \in \mathbb{H}^3 \mid -\frac{1}{4} \leq x \leq 0, \; y = \frac{1+n^2d+4x}{4n\sqrt{d}} \right\}\]
\[F_U = \left\{(x,y,t) \in \mathbb{H}^3 \mid \frac{1}{4} \leq x \leq \frac{1}{2}, \; y = \frac{n^2d+1-4x}{4n\sqrt{d}} \right\}\]
\[F_{U^{-1}} = \left\{(x,y,t) \in \mathbb{H}^3 \mid 0 \leq x \leq \frac{1}{4}, \; y = \frac{n^2d+1-4x}{4n\sqrt{d}} \right\}\]
and edges

\begin{tabular}{lll}
$\bullet \; e_{(1),\pi} = F_S \cap F_{S^{-1}}$&$\bullet \; e_{(3,1),\theta} = F_A \cap F_T$&$\bullet \; e_{(4),\pi} = F_T \cap F_{T^{-1}}$\\
$\bullet \; e_{(2,1),\pi/3} = F_A \cap F_S$&$\bullet \; e_{(3,2),\theta} = F_{A^{-1}} \cap F_U$&$\bullet \; e_{(5),\pi} = F_U \cap F_{U^{-1}}$\\
$\bullet \; e_{(2,2),\pi/3} = F_{A^{-1}} \cap F_{S^{-1}}$&$\bullet \; e_{(3,3),\phi} = F_{T^{-1}} \cap F_{U^{-1}}$\\
\end{tabular}

\begin{figure}[htp]
\begin{center}
\input{subcase222a.pstex_t}
\end{center}
\caption[A fundamental polyhedron in Subcase 2.2.2 of Theorem 3.16]{The polygon $P$ viewed from the point at $\infty$ above.  Solid lines denote edges, solid dots denote edges parallel to the $t$-axis, dashed lines are on the boundary of $\mathbb{H}^3$, and dotted lines denote the $x$- and $y$-axes in the boundary of $\mathbb{H}^3$.  The sides $F_S$ and $F_{S^{-1}}$ are on the unit sphere centered at the origin; all other sides are perpendicular to the boundary of $\mathbb{H}^3$.}
\label{subcase222a}
\end{figure}

Then $X(F_X) = F_{X^{-1}}$ for $X \in \{A,S,T,U\}$ are side pairing transformations with the following five cycle transformations:
\begin{enumerate}
\item $e_{(1),\pi} \stackrel{S}{\longrightarrow} \circlearrowleft^2_\pi$
\item $e_{(2,1),\pi/3} \stackrel{A}{\longrightarrow} e_{(2,2),\pi/3} \stackrel{S^{-1}}{\longrightarrow}  \circlearrowleft^3_{2\pi/3}$
\item $e_{(3,1),\theta} \stackrel{A}{\longrightarrow} e_{(3,2),\theta} \stackrel{U}{\longrightarrow} e_{(3,3),\phi} \stackrel{T^{-1}}{\longrightarrow} \circlearrowleft^1_{2\pi}$
\item $e_{(4),\pi} \stackrel{T}{\longrightarrow} \circlearrowleft^2_\pi$
\item $e_{(5),\pi} \stackrel{U}{\longrightarrow} \circlearrowleft^2_\pi$
\end{enumerate}
Therefore, by Poincar\'{e}'s Polyhedron Theorem, $\Gamma$ has the presentation
\[ \Gamma = \langle A,S,T,U \mid S^2 = (S^{-1}A)^3 = T^{-1}UA = T^2 = U^2 = 1 \rangle, \]
and $P$ is the fundamental polyhedron for $\Gamma$. Clearly $P$ has infinite volume, so $\Gamma$ is not arithmetic if $n>1$ or $d>11$.

The case of $n=1$ and $d=3, 7$, or $11$ is similar (i.e., $\Gamma$ does not have finite covolume and hence is not arithmetic); see \cite{Callahan09} for details.

\textbf{Case 3.} Suppose $\sigma = e^{-\pi i/4}$ and $d=1$. Then $\tr B \tr AB - \tr^2 B = \lambda e^{-\pi i/4} \in O_1$.  Thus, $\lambda e^{-\pi i/4} = m + n e^{-\pi i/2}$ for some integers $m$ and $n$, so $\lambda = m e^{\pi i/4} + ne^{-\pi i/4}$.  Replacing $B$ with $BA^{-n}$, assume $\lambda = m e^{\pi i/4}$.  Then $\Gamma = G_{\theta, k}$ with $(\theta, k) = (\frac{\pi}{4}, \frac{m}{2})$, so, by Proposition  \ref{losid}, $m = 1$, 2, or 3, in which case $\Gamma$ is $\mathrm{PGL}_2(O_1)$, a subgroup of index 8 in $\mathrm{PGL}_2(O_1)$, or a subgroup of index 6 in $\mathrm{PGL}_2(O_1)$ respectively.

The case of $\sigma = e^{\pi i/4}$ and $d=1$ is similar; see \cite{Callahan09} for details.

\textbf{Case 4.} Suppose $\sigma = e^{-\pi i/6}$ and $d=3$.  Then
\[(\tr B \tr AB - \tr^2 B)(\tr[A,B] - 2) = \lambda e^{-\pi/2} \in O_3. \] Hence, $\lambda e^{-\pi/2} = m + n e^{-2\pi/3}$ for some integers $m$ and $n$, so $\lambda = mi + ne^{-\pi i/6}$.  Replacing $B$ with $BA^{-n}$, assume $\lambda = mi$.

If $m$ is even, then replace $B$ by $BA^{m/2}$ and assume
\[\lambda = mi + \frac{m}{2}e^{-\pi i/6} = \frac{m\sqrt{3}}{2}e^{\pi i/3}.\]
Hence, $\Gamma = G_{\theta, k}$ with $(\theta, k) = (\frac{\pi}{3}, \frac{m\sqrt{3}}{4})$, so $m$ can be any even integer by Proposition  \ref{losid}; $\Gamma$ is a $\Z_2$-extension of the figure-eight knot group by an involution that conjugates one parabolic generator to the other if $m \equiv 0 \mod 4$ and by an involution that conjugates each parabolic generator to its own inverse if $m \equiv 2 \mod 4$.

Now suppose $m$ is odd.  Let
\begin{eqnarray*}
C &=& BA^{-1}B^{-1} = \left(\begin{array}{cc} 1&0 \\ e^{-\pi i/3}&1 \\ \end{array}\right) \text{ and }\\
T &=& C^{-1}ACA^{-2}CAC^{-1} = \left(\begin{array}{cc} 1&-2\sqrt{-3} \\ 0&1 \\ \end{array}\right).
\end{eqnarray*}

If $m \equiv 1 \mod 4$, then let
\[D = BA^{\frac{m+3}{2}}T^{\frac{m-1}{4}} = \left(\begin{array}{cc} 0&-e^{\pi i/6} \\ e^{-\pi i/6}&\sqrt{3} \\ \end{array}\right),\]
so $B \in \langle A, D \rangle$ since $T \in \langle A, C \rangle$ and $C = DAD^{-1}$.  Since $D$ is elliptic of order 6, we conclude that $\Gamma = \langle A, D \rangle$ is $\mathrm{PGL}_2 (O_3)$ by Theorem \ref{conder}.  Thus, $\mathrm{PGL}_2 (O_3)$ is an arithmetic J\o{}rgensen group of parabolic type but is not identified in Proposition  \ref{losid}.  Hence, $\mathrm{PGL}_2 (O_3)$ is not of the form $G_{\theta, k}$.

The cases when $m \equiv 3 \mod 4$ and when $\sigma = e^{\pi i/6}$ and $d=3$ are similar; see \cite{Callahan09} for details.

\textbf{Case 5.}  Suppose $\sigma = e^{-\pi i/3}$ and $d=3$.  Then
\[(\tr B \tr AB - \tr^2 B)(\tr[A,B] - 2) = -\lambda \in O_3. \] Hence, $\lambda = m+ne^{-\pi i/3}$ for some integers $m$ and $n$.  Replacing $B$ with $BA^{-n}$, assume $\lambda = m$.

If $m$ is even, then replace $B$ by $BA^{-m/2}$ and assume
\[\lambda = m - \frac{m}{2}e^{-\pi i/3} = \frac{m\sqrt{3}}{2}e^{\pi i/6}.\]
Hence, $\Gamma = G_{\theta, k}$ with $(\theta, k) = (\frac{\pi}{6}, \frac{m\sqrt{3}}{4})$, so $m$ can be any even integer by Proposition  \ref{losid}; $\Gamma$ is the figure-eight knot group if $m \equiv 2 \mod 4$ and is a $\Z_2$-extension of the figure-eight knot group by an involution that conjugates one parabolic generator to the other if $m \equiv 0 \mod 4$ (recall that this $\Z_2$-extension of the figure-eight knot group is not conjugate to the one found in Case 4 by comparing trace fields).

Now suppose $m$ is odd.  Let
\begin{eqnarray*}
C &=& BAB^{-1} = \left(\begin{array}{cc} 1&0 \\ e^{\pi i/3}&1 \\ \end{array}\right) \text{ and }\\
T &=& C^{-1}ACA^{-2}CAC^{-1} = \left(\begin{array}{cc} 1&2\sqrt{-3} \\ 0&1 \\ \end{array}\right).
\end{eqnarray*}

If $m \equiv 1 \mod 4$, then let
\[D = A^{-1}CAC^{-1}B^{-1}A^{-1}CAC^{-1}A^{\frac{m-1}{2}}T^{\frac{m-5}{4}} = \left(\begin{array}{cc} 0&-e^{\pi i/3} \\ e^{-\pi i/3}&-1 \\ \end{array}\right),\]
so $B \in \langle A, D \rangle$ since $T \in \langle A, C \rangle$ and $C = DAD^{-1}$.  Since $D$ is elliptic of order 3, we conclude that $\Gamma = \langle A, D \rangle$ is $\mathrm{PSL}_2 (O_3)$ by Theorem \ref{conder}.  Thus, $\mathrm{PSL}_2 (O_3)$ is an arithmetic J\o{}rgensen group of parabolic type but is not identified in Proposition  \ref{losid}.  Hence, $\mathrm{PSL}_2 (O_3)$ is not of the form $G_{\theta, k}$.

The cases when $m \equiv 3 \mod 4$ and when $\sigma = e^{\pi i/3}$ and $d=3$ are similar; see \cite{Callahan09} for details.
\end{proof}

\section{Arithmetic J\o{}rgensen Groups of Elliptic Type} \label{elliptic}

Before discussing arithmetic J\o{}rgensen groups of elliptic type, we recall several more facts about two-generator Kleinian groups and arithmeticity.  The first, Theorem 3.6.2 in \cite{MaclachlanReid03}, will enable us to determine invariant quaternion algebras.

\begin{thm} \label{quatalg}
If $\langle X, Y \rangle$ is an irreducible subgroup of a non-elementary group $\Gamma$ in $\pslc$ such that $X$ is not parabolic and neither $X$ nor $Y$ is elliptic of order two, then the invariant quaternion algebra is \[ A\Gamma = \left( \dfrac{\tr^2 X (\tr^2 X - 4), \; \tr^2 X \tr^2 Y (\tr [X,Y] - 2)}{k\Gamma} \right).\]
\end{thm}

We next state Theorem 8.3.2 of \cite{MaclachlanReid03}, together with Lemma 4.1 in \cite{GehringMaclachlanMartinReid97}, which will enable us to determine when a finite-covolume Kleinian group is arithmetic.

\begin{thm} \label{idthm}
A finite-covolume Kleinian group $\Gamma$ is arithmetic if and only if the following three conditions hold.
\begin{enumerate}
\item $k\Gamma$ is a number field with exactly one complex place.
\item $\tr \Gamma$ consists of algebraic integers.
\item $A\Gamma = \left( \dfrac{a, b}{k\Gamma} \right)$ is ramified at all real places of $k\Gamma$, i.e., $\tau(a)$ and $\tau(b)$ are negative for all real embeddings $\tau$ of $k\Gamma$.
\end{enumerate}
\end{thm}

Finally, turning our attention to Fuchsian subgroups of arithmetic Kleinian groups, we note the following combination of Theorem 9.5.2 and Corollary 9.5.3 in \cite{MaclachlanReid03} (cf. Theorems 4 and 5 of \cite{MaclachlanReid87}).

\begin{thm} \label{fuch}
If $F$ is a non-elementary Fuchsian subgroup of an arithmetic Klein-ian group $\Gamma$, then $F$ is a subgroup of an arithmetic Fuchsian group $G$ such that \[kG = k\Gamma \cap \R, \; [k\Gamma:kG] = 2, \text{ and } A\Gamma \cong AG \otimes_{kG} k\Gamma.\]
\end{thm}

We now establish the following characterization of arithmetic J\o{}rgensen groups of elliptic type, which completes the proof of Theorem \ref{parabcase} as a corollary.

\begin{thm} \label{ellipcase}
A finite-covolume Kleinian group $\Gamma$ is an arithmetic J\o{}rgensen group of elliptic type if and only if $\Gamma = \langle A,B \rangle$ such that the following six conditions hold.
\begin{enumerate}
\item $A$ is elliptic of order $n \in \{7, 8, 9, 10, 11, 12, 14, 16, 18, 24, 30\}$.
\item $B$ is loxodromic or hyperbolic with $\tr^2 B > \dfrac{2}{1-\cos\frac{2\pi}{n}} > 4$.
\item $\tr [A,B] = 2\cos\frac{2\pi}{n} + 1$.
\item $k\Gamma$ is a quadratic imaginary extension of $\Q(\cos\frac{2\pi}{n})$ such that \[ -1 < \tau\left(\cos\frac{2\pi}{n}\right) < \frac{1}{2} \text{ and } 0 < \tau(\tr^2 B) < \dfrac{2}{1-\tau(\cos\frac{2\pi}{n})} < 4 \] for each of the $2[\Q(\cos\frac{2\pi}{n}):\Q]-2$ real embeddings $\tau$ of $k\Gamma$.
\item $A\Gamma \cong \left( \dfrac{-1, \; 2\cos\frac{2\pi}{n}-1}{k\Gamma} \right) \cong \left( \dfrac{-1, \; 2(\cos\frac{4\pi}{n} + \cos\frac{2\pi}{n}) \tr^2 B} {k\Gamma} \right)$.
\item $\tr B$ and $\tr AB$ are algebraic integers.
\end{enumerate}
\end{thm}

\begin{proof}
We first prove necessity.  Then $\Gamma = \langle A,B \rangle$ is an arithmetic Kleinian group with $J(A,B)=1$ and $A$ elliptic. By Theorem \ref{kiikka}, $A$ has order $n \geq 7$ and $\tr ABAB^{-1} = 1$.  Thus, $\tr^2 A = 2\cos\frac{2\pi}{n} + 2$, and standard trace relations (3.14 of \cite{MaclachlanReid03}, for instance) yield
\[ 1 = \tr ABAB^{-1} = (\tr A)(\tr BAB^{-1}) - \tr A(BAB^{-1})^{-1} = \tr^2 A - \tr [A,B], \]
so $\tr [A,B] = 2\cos\frac{2\pi}{n} + 1$, which establishes Condition (3) of the theorem.

Also by Theorem \ref{kiikka}, $\Delta = \langle A, BAB^{-1} \rangle$ is a non-elementary subgroup of $\Gamma$ with $J(A, BAB^{-1}) = 1$, so $\tr [A, BAB^{-1}] = 2\cos\frac{2\pi}{n} + 1$ as above.  Since \[\tr A = \tr BAB^{-1} = \pm 2\cos\frac{\pi}{n} \text{ and } \tr ABAB^{-1} = 1,\] Lemma \ref{tracefields} yields $k\Delta = \Q(\cos\frac{2\pi}{n})$ and $\Q(\tr\Delta) = \Q(\cos\frac{\pi}{n}) \subset \R$, so $\Delta$ is a Fuchsian subgroup of $\Gamma$ (cf. Corollary 3.2.5 of \cite{MaclachlanReid03}).  Therefore, by Theorem \ref{fuchjorg}, $\Delta$ must be a $(2,3,q)$-triangle group, which has trace field $\Q(\cos\frac{\pi}{q})$ (see, for instance, Section 4.9 of \cite{MaclachlanReid03}), forcing $n=q$.

Theorem \ref{fuch} further implies that $\Delta$ is a subgroup of an arithmetic Fuchsian group $G$.  But, by Theorem 3B of \cite{Greenberg63}, $(2,3,n)$-triangle groups cannot be subgroups of strictly lager Fuchsian groups.  Therefore, $\Delta = G$, and so $\Delta$ is an arithmetic $(2,3,n)$-triangle group.  Using the enumeration of all arithmetic triangle groups in \cite{Takeuchi77} and noting that $n \neq \infty$ since $A$ is elliptic of order $n$, the $(2,3,n)$-triangle group is arithmetic if and only if \[n \in \{7, 8, 9, 10, 11, 12, 14, 16, 18, 24, 30\}\] (cf. Section 13.3 in \cite{MaclachlanReid03}), thereby establishing Condition (1).

By Theorems \ref{idthm} and \ref{fuch}, $\Q(\cos\frac{2\pi}{n}) = k\Gamma \cap \R$, $[k\Gamma:\Q(\cos\frac{2\pi}{n})] = 2$, and $k\Gamma$ has exactly one complex place, so $k\Gamma$ must be a quadratic imaginary extension of $\Q(\cos\frac{2\pi}{n})$.  Applying Theorem \ref{quatalg} to $\langle A, BAB^{-1} \rangle$ yields
\begin{eqnarray*}
A\Gamma &=& \left( \dfrac{(2\cos\frac{2\pi}{n} + 2)(2\cos\frac{2\pi}{n} - 2), \; (2\cos\frac{2\pi}{n} + 2)^2(2\cos\frac{2\pi}{n}-1)}{k\Gamma} \right) \\
&\cong& \left( \dfrac{-1, \; 2\cos\frac{2\pi}{n}-1}{k\Gamma} \right) \\
\end{eqnarray*}
where we use Lemma 2.1.2 in \cite{MaclachlanReid03} to remove squares of elements in $k\Gamma^*$. Similarly, applying Theorem \ref{quatalg} to $\langle A, B \rangle$ yields
\begin{eqnarray*}
A\Gamma &=& \left( \dfrac{(2\cos\frac{2\pi}{n} + 2)(2\cos\frac{2\pi}{n} - 2), \; (2\cos\frac{2\pi}{n} + 2)\tr^2 B(2\cos\frac{2\pi}{n}-1)}{k\Gamma} \right) \\
&\cong& \left( \dfrac{-1, \; 2(\cos\frac{4\pi}{n} + \cos\frac{2\pi}{n})\tr^2 B}{k\Gamma} \right) \\
\end{eqnarray*}
again using Lemma 2.1.2 in \cite{MaclachlanReid03} to remove squares of elements in $k\Gamma^*$.  This establishes Condition (5).

Clearly $k\Gamma$ has $2[\Q(\cos\frac{2\pi}{n}):\Q]-2$ real places, each of which corresponds to a real embedding.  Let $\tau: k\Gamma \rightarrow \R$ be one such real embedding.  Then, since
\[ A\Gamma \cong \left( \dfrac{(2\cos\frac{2\pi}{n} + 2)(2\cos\frac{2\pi}{n} - 2), \; 2\cos\frac{2\pi}{n}-1}{k\Gamma} \right) \] is ramified at all real places of $k\Gamma$ by Theorem \ref{idthm}, we have
\[ \tau\left(\cos\frac{2\pi}{n}\right) - 1 < 0 < \tau\left(\cos\frac{2\pi}{n}\right) + 1 \text{ and } 2\tau\left(\cos\frac{2\pi}{n}\right) - 1 < 0 \]
which yields $-1 < \tau(\cos\frac{2\pi}{n}) < \frac{1}{2}$.  Similarly, for
\[ A\Gamma \cong \left( \dfrac{-1, (2\cos\frac{2\pi}{n} + 2)\tr^2 B(2\cos\frac{2\pi}{n}-1)}{k\Gamma} \right) \] to be ramified at all real places of $k\Gamma$, we must have $\tau(\tr^2 B) > 0$ since
\[ \tau\left(2\cos\frac{2\pi}{n} - 1\right) < 0 < \tau\left(2\cos\frac{2\pi}{n} + 2\right) \] as above.

By Lemma \ref{tracefields}, $k\Gamma = \mathbb{Q}(\tr^2 A, \tr^2 B, \tr A \tr B \tr AB)$.  Following Section 4 of \cite{MaclachlanMartin99}, the standard trace relation (3.15 in \cite{MaclachlanReid03}, for instance)
\[ \tr [A,B] = \tr^2 A + \tr^2 B + \tr^2 AB - \tr A \tr B \tr AB - 2 \]
implies that $\tr A \tr B \tr AB$ satisfies the quadratic equation
\[ x^2 - (\tr^2 A )(\tr^2 B)x - (\tr^2 A )(\tr^2 B)(\tr [A,B] - \tr^2 A - \tr^2 B + 2) = 0, \]
so \[ k\Gamma = \mathbb{Q}(\tr^2 A, \tr^2 B, \tr [A,B])(\sqrt{\delta}) \]
where
\begin{eqnarray*}
\delta &=& (\tr^2 A )(\tr^2 B)[(\tr^2 A - 4)(\tr^2 B - 4)+4(\tr [A,B]-2)] \\
&=& \left(2\cos\frac{2\pi}{n} + 2\right)(\tr^2 B)\left[\left(2\cos\frac{2\pi}{n} - 2\right)(\tr^2 B - 4)+4\left(2\cos\frac{2\pi}{n} - 1\right)\right]. \\
\end{eqnarray*}
As $\tau$ is a real embedding of $k\Gamma$, we must have $\tau(\delta) > 0$. Having already shown that
\[-1 < \tau\left(\cos\frac{2\pi}{n}\right) < \frac{1}{2} \text{ and } \tau(\tr^2 B) > 0, \]
we conclude that
\[ \tau(\tr^2 B) < \dfrac{2}{1-\tau(\cos\frac{2\pi}{n})} < 4.\]  Condition (4) is now established.

Suppose $\tr B$ is real. Then $\tr^2 B \in \R \cap k\Gamma = \Q(\cos\frac{2\pi}{n})$, so $k\Gamma = \Q(\cos\frac{2\pi}{n})(\sqrt{\delta})$ as above.  If $\tr^2 B = 0$, then $\delta = 0$, which implies $k\Gamma \subset \R$, a contradiction, so $\tr^2 B > 0$. Since $k\Gamma$ has exactly one complex place, we conclude that $\delta < 0$, and so
\[ (\tr^2 A - 4)(\tr^2 B - 4)+4(\tr [A,B]-2) = (2\cos\frac{2\pi}{n} - 2)(\tr^2 B - 4)+4(2\cos\frac{2\pi}{n} - 1) < 0 \]
since $\tr^2 A > 0$ as well.  Thus, \[ \tr^2 B > \dfrac{2}{1-\cos\frac{2\pi}{n}} > 4 \]
since $n > 6$, so $B$ cannot be elliptic or parabolic.  Condition (2) is now established, and Condition (6) follows directly from Theorem \ref{idthm}, so necessity has been proved.

We now verify that sufficiency follows from construction.  Suppose $\Gamma = \langle A,B \rangle$ is a finite-volume Kleinian group such that the six conditions of the theorem hold. Then Condition (4) ensures that $k\Gamma$ is a number field with exactly one complex place, Conditions (1) and (6) imply that $\tr \Gamma$ consists of algebraic integers by Lemma \ref{tracefields}, and Conditions (4) and (5) guarantee that $A\Gamma$ is ramified at all real places of $k\Gamma$, so $\Gamma$ is arithmetic by Theorem \ref{idthm}.  Finally, Conditions (1) and (3) yield
\begin{eqnarray*}
J(A,B) &=& |\tr^2 A - 4| + |\tr [A,B] - 2| \\
       &=& \left|2\cos\frac{2\pi}{n} - 2\right| + \left|2\cos\frac{2\pi}{n} - 1\right| \\
       &=& 2 - 2\cos\frac{2\pi}{n} + 2\cos\frac{2\pi}{n} - 1 \\
       &=& 1 \\
\end{eqnarray*}
since $n>6$. Thus, $\Gamma$ is an arithmetic J\o{}rgensen group of elliptic type.
\end{proof}

\begin{cor}
No arithmetic J\o{}rgensen group of elliptic type is commensurable with any non-compact arithmetic Kleinian group.
\end{cor}
 
\begin{proof}
By Theorem \ref{bianchi}, the invariant trace field of a non-cocompact arithmetic Kleinian group has degree 2 over $\Q$, whereas by Condition (4) of Theorem \ref{ellipcase}, the invariant trace field of an arithmetic J\o{}rgensen group of elliptic type has degree $2[\Q(\cos\frac{2\pi}{n}):\Q]$ over $\Q$, where $n$ is as in Condition (1) of the theorem. Since the invariant trace field is an invariant of the commensurability class of a finitely generated non-elementary Kleinian group (see, for instance, Section 3.3 of \cite{MaclachlanReid03}), the result follows.
\end{proof}

Thus, the non-cocompact arithmetic J\o{}rgensen groups are precisely the arithmetic J\o{}rgensen groups of parabolic type identified in Theorem \ref{parabtype}, thereby completing the proof of Theorem \ref{parabcase}.  We also note that Conditions (1) and (3) of Theorem \ref{ellipcase} show that arithmetic J\o{}rgensen groups of elliptic type are Kleinian groups with real parameters, which are discussed in \cite{GehringGehringMartin01}.

\section{Bounds on J\o{}rgensen Number} \label{bounds}

We begin with a simple bound on J\o{}rgensen numbers of Kleinian groups containing parabolic elements.

\begin{pro} \label{genbds}
Let $\Gamma$ be a Kleinian group that contains the parabolic element $A = \left(\begin{array}{cc} 1&1 \\ 0&1 \\ \end{array}\right)$.  Then 
\[ 1 \leq \widetilde{J}(\Gamma) \leq \inf \left\{|c|^2 : 
T = \left(\begin{array}{cc} a&b \\ c&d \\ \end{array}\right) \in \Gamma, \; c \neq 0, \text{ and } 
\langle A, T \rangle \text{ is non-elementary} \right\}. \] 
If $\Gamma$ is non-elementary and can be generated by $A$ and another element, then
\[ 1 \leq J(\Gamma) \leq \inf \left\{|c|^2 : 
T = \left(\begin{array}{cc} a&b \\ c&d \\ \end{array}\right) \in \Gamma, \; c \neq 0, \text{ and } 
\langle A, T \rangle = \Gamma \right\}. \]
\end{pro}

\begin{proof}
As already noted, $1 \leq \widetilde{J}(\Gamma) \leq J(\Gamma)$. Let 
$T = \left(\begin{array}{cc} a&b \\ c&d \\ \end{array}\right) \in \Gamma$ with $c \neq 0.$ If $\langle A, T
\rangle$ is non-elementary, then $\widetilde{J}(\Gamma) \leq J(A,T) = |c|^2$.  If $\langle A, T \rangle = \Gamma$, then $J(\Gamma) \leq J(A,T) = |c|^2$.  Since $T$ was arbitrary, the bounds follow.
\end{proof}

Now suppose that $\mathbb{H}^3/\Gamma$ is a hyperbolic two-bridge knot complement $S^3 \smallsetminus K$. Following Section 3 of \cite{Riley72} and Section 4.5 of \cite{MaclachlanReid03}, $K$ is determined by a pair of relatively prime odd integers $(p/q)$ with $0 < q < p$, and $\Gamma = \langle A, B \; | \; AW = WB \rangle$ where
\[ A = \left(\begin{array}{cc} 1&1 \\ 0&1 \\ \end{array}\right), \; B = \left(\begin{array}{rr} 1&0 \\ -z&1 \\\end{array}\right), \; W = B^{e_1}A^{e_2} \cdots B^{e_{p-2}}A^{e_{p-1}},\] and $e_i = (-1)^{\left \lfloor \frac{iq}{p} \right \rfloor}$.

The relation $AW = WB$ forces 
$W = \left(\begin{array}{cc} \ast & 1/\sqrt{z} \\ -\sqrt{z}&0 \\\end{array}\right)$ where $z$ satisfies the integral monic polynomial equation 
\[ d_n = 1 + \left(\widehat{\sum_{i_1\text{ odd}}}e_{i_1}e_{i_2} \right) z + \cdots + (e_1e_2\cdots e_{2n})z^n = 0\] where $n = \dfrac{p-1}{2}$ and $\widehat{\sum}$ denotes summation over $i_1 < i_2 < \cdots < i_k$ with alternating parity.  Then $k\Gamma = \Q(\tr \Gamma) = \Q(z)$, and since two-bridge knot groups are not free, a well known result establishes the bound $|z| < 4$ (see, for instance, Theorem B of \cite{Riley72}, which is attributed to J.\ Brenner).

\begin{cor}
Suppose $M = \mathbb{H}^3/\Gamma$ is the complement in $S^3$ of a hyperbolic two-bridge knot $K$ with all notation as above, and let $\mathcal{C}$ be the single cusp in $M$. If $K$ is the figure-eight knot, then \[ \widetilde{J}(\Gamma) = J(\Gamma) = 1 = w(M, \mathcal{C}).\] Otherwise, \[ 1 < \widetilde{J}(\Gamma) \leq J(\Gamma) \leq |z| < 4 \text{ and } 1 < w(M, \mathcal{C}) \leq |\sqrt{z}| < 2. \]
\end{cor}

\begin{proof}
Results for the figure-eight knot have already been established in Section \ref{manifolds}, so let $K$ be any other hyperbolic two-bridge knot $(p/q)$ with notation as above.  The relation $AW = WB$ yields $B = W^{-1}AW$, so $\langle A, W \rangle = \Gamma$, and $J(A,W) = |z|$.  Hence, \[ 1 < \widetilde{J}(\Gamma) \leq J(\Gamma) \leq |z| < 4. \]

The bounds on $w(M, \mathcal{C})$ follow from Theorem \ref{adams} and Lemma \ref{waistlem} applied to $W$.
\end{proof}

Similarly, suppose that $\mathbb{H}^3/\Gamma$ is a hyperbolic two-bridge link complement $S^3 \smallsetminus L$ where $L$ has two components.
Again following Section 4.5 of \cite{MaclachlanReid03}, $L$ is determined by a pair of relatively prime integers $(p/q)$ with $0 < q < p$, $p = 2n$ even, and $\Gamma = \langle A, B \; | \; AW = WA \rangle$ where
\[ A = \left(\begin{array}{cc} 1&1 \\ 0&1 \\ \end{array}\right), \; B = \left(\begin{array}{rr} 1&0 \\ -z&1 \\\end{array}\right), \; W = B^{e_1}A^{e_2} \cdots B^{e_{p-1}},\] and $e_i = (-1)^{\left \lfloor \frac{iq}{p} \right \rfloor}$.

The relation $AW = WA$ forces 
$z$ to satisfy the integral monic polynomial equation \[ c_n = \left(\sum_{i=1}^n e_{2i-1}\right) z + \left( \widehat{\sum_{i_1\text{ odd}}}e_{i_1}e_{i_2}e_{i_3} \right) z^2 + \cdots + (e_1e_2\cdots e_{2n-1})z^n = 0 \] where $\widehat{\sum}$ again denotes summation over $i_1 < i_2 < \cdots < i_k$ with alternating parity.  As before, $k\Gamma = \Q(\tr \Gamma) = \Q(z)$, and $|z| < 4$.

\begin{cor}
Suppose $M = \mathbb{H}^3/\Gamma$ is the complement in $S^3$ of a hyperbolic two-bridge link $L$ with all notation as above, and let $\mathcal{C}$ be the cusp in $M$ whose lift contains a horoball based at $\infty$. Then
\[ 1 < \widetilde{J}(\Gamma) \leq J(\Gamma) \leq |z|^2 < 16 \text{ and } 1 < w(M, \mathcal{C}) \leq |z| < 4. \]
\end{cor}

\begin{proof}
Since $L$ has two components, it is not the figure-eight knot, so the bounds follow as before with Proposition \ref{genbds} and Lemma \ref{waistlem} applied to $B$ this time.
\end{proof}

\section{Computations of J\o{}rgensen Numbers} \label{comp}

To compute the J\o{}rgensen numbers of several non-cocompact Kleinian groups, we first collect some preliminary results that will be useful in what follows. We begin with Lemma 7.1, combined the comments and definitions that precede it, in \cite{GehringMaclachlanMartinReid97}.

\begin{lem} \label{GMMR}
Let $\Gamma$ be a finite-covolume Kleinian group whose traces lie in $R$, the ring of integers in $\mathbb{Q}(\tr\Gamma)$.  If $\langle X, Y \rangle$ is a non-elementary subgroup of $\Gamma$, then $\mathcal{O} = R[1, X, Y, XY]$ is an order in the quaternion algebra
\[ A_0\Gamma = \left\{ \sum a_i \gamma_i : a_i \in \mathbb{Q}(\tr\Gamma), \gamma_i \in \Gamma \right\} \]
over $\mathbb{Q}(\tr\Gamma)$. Its discriminant $d(\mathcal{O})$ is the ideal $\langle \tr [X, Y] - 2 \rangle$ in $R$.
\end{lem}

The second is Theorem 6.3.4 in \cite{MaclachlanReid03}.

\begin{thm} \label{MR}
Let $\mathcal{O}_1$ and $\mathcal{O}_2$ be orders in a quaternion algebra over a number field.  If $\mathcal{O}_1 \subset \mathcal{O}_2$, then $d(\mathcal{O}_2) \, | \, d(\mathcal{O}_1)$, and $\mathcal{O}_1 = \mathcal{O}_2$ if and only if $d(\mathcal{O}_1) = d(\mathcal{O}_2)$.
\end{thm}

Our application is the following.

\begin{cor} \label{unitmult}
Let $\Gamma$ be a finite-covolume Kleinian group whose traces lie in $R$, the ring of integers in $\mathbb{Q}(\tr\Gamma)$.  If $\langle A, B \rangle = \Gamma = \langle X, Y \rangle$, then $\tr [X,Y] - 2$ is a unit multiple of $\tr [A,B] - 2$ in $R$.
\end{cor}

\begin{proof}
By Lemma \ref{GMMR}, $\mathcal{O}_1 = R[1, A, B, AB]$ and $\mathcal{O}_2 = R[1, X, Y, XY]$ are orders in $A_0\Gamma$.  Furthermore, $d(\mathcal{O}_1) = \langle 2- \tr[A, B] \rangle$ and $d(\mathcal{O}_2) = \langle 2 - \tr[X, Y] \rangle$ are ideals in $R$.  The Cayley-Hamilton Theorem yields the identity
\[ X + X^{-1} = \tr X \cdot 1, \]
which implies $A^{-1}, B^{-1} \in \mathcal{O}_1$ and $X^{-1}, Y^{-1} \in \mathcal{O}_2$.  Since $\mathcal{O}_1$ and $\mathcal{O}_2$ are ideals that are also rings with 1, we have
\[\Gamma = \langle A, B \rangle \subset \mathcal{O}_1 \text{ and } \Gamma = \langle X, Y \rangle \subset \mathcal{O}_2.\]
Let \[ R\Gamma = \left\{ \sum a_i \gamma_i : a_i \in R, \gamma_i \in \Gamma \right\}. \] Clearly $R \subset \mathcal{O}_1,\mathcal{O}_2$, so $R\Gamma \subseteq \mathcal{O}_1,\mathcal{O}_2$, and $\mathcal{O}_1,\mathcal{O}_2 \subseteq R\Gamma$ by definition.  Therefore,
\[\mathcal{O}_1 = R\Gamma = \mathcal{O}_2.\]
Hence, $d(\mathcal{O}_1) = d(\mathcal{O}_2)$ by Theorem \ref{MR}, and the result follows.
\end{proof}

The following proposition yields as a corollary the computation of the J\o{}rgensen number of every two-generator arithmetic Kleinian group with one generator parabolic and the other parabolic or elliptic.

\begin{pro} \label{arithcomp}
If $\Gamma = \langle A,B \rangle$ is an arithmetic Kleinian group with $A = \left(\begin{array}{cc} 1&1 \\ 0&1 \\ \end{array}\right)$ and $B$ parabolic or elliptic, then $J(\Gamma) = J(A,B)$.
\end{pro}

\begin{proof}
Suppose $\Gamma = \langle X,Y \rangle$.  Then, by Lemma \ref{tracefields}, Theorem \ref{bianchi}, and Corollary \ref{unitmult}, $\tr [X,Y] - 2$ is a unit multiple of $\tr [A,B] - 2$ in $O_d$, all of whose units have norm one.  Hence, \[ J(A,B) = |\tr [A,B] - 2| = |\tr [X,Y] - 2| \leq J(X,Y).\] The result now follows.
\end{proof}

\begin{cor}
The J\o{}rgensen numbers of every two-generator arithmetic Kleinian group with one generator parabolic and the other parabolic or elliptic are as follows.
\begin{itemize}
\item If $\Gamma$ is the figure-eight knot group, then $J(\Gamma) = 1$.
\item If $\Gamma$ is the Whitehead link group, then $J(\Gamma) = 2$.
\item If $\Gamma$ is the $6_2^2$ link group, then $J(\Gamma) = 3$.
\item If $\Gamma$ is the $6_3^2$ link group, then $J(\Gamma) = 2$.
\item If $\Gamma$ is a $\mathbb{Z}_2$-extension of the figure eight knot group, then $J(\Gamma) = 1$.
\item If $\Gamma$ is a $\mathbb{Z}_2$-extension of the Whitehead link group, then $J(\Gamma) = \sqrt{2}$.
\item If $\Gamma$ is a $\mathbb{Z}_2$-extension of the $6^2_2$ link group, then $J(\Gamma) = \sqrt{3}$.
\item If $\Gamma$ is a $\mathbb{Z}_2$-extension of the $6^2_3$ link group, then $J(\Gamma) = \sqrt{2}$.
\item If $[\mathrm{PSL}_2 (O_1) : \Gamma] = 8$, then $J(\Gamma) = 2$.
\item If $\Gamma = \mathrm{PSL}_2 (O_3)$, then $J(\Gamma) = 1$.
\item If $[\mathrm{PSL}_2 (O_7) : \Gamma] = 2$, then $J(\Gamma) = 2$.
\item If $\Gamma = \mathrm{PGL}_2 (O_1)$, then $J(\Gamma) = 1$.
\item If $[\mathrm{PSL}_2 (O_2) : \Gamma \cap \mathrm{PSL}_2 (O_2)] = 24$, then $J(\Gamma) = 3$.
\item If $[\mathrm{PSL}_2 (O_3) : \Gamma \cap \mathrm{PSL}_2 (O_3)] = 30$, then $J(\Gamma) = 2$.
\item If $\Gamma = \mathrm{PGL}_2 (O_3)$, then $J(\Gamma) = 1$.
\item If $[\mathrm{PSL}_2 (O_{15}) : \Gamma \cap \mathrm{PSL}_2 (O_{15})] = 6$, then $J(\Gamma) = 2$.
\end{itemize}
\end{cor}

\begin{proof}
As established in \cite{GehringMaclachlanMartin98} and \cite{Conderetal02} (cf. Theorem \ref{conder}), every arithmetic Kleinian group $\Gamma = \langle A,B \rangle$ with
\[ A = \left(\begin{array}{cc} 1&1 \\ 0&1 \\ \end{array}\right) \text{ and } B = \left(\begin{array}{cc} a&b \\ c&d \\ \end{array}\right), \] where $B$ is parabolic or elliptic, is listed above and $c$ is as follows.
\begin{itemize}
\item If $\Gamma$ is the figure-eight knot group, then $c = \dfrac{1+i\sqrt{3}}{2}$.
\item If $\Gamma$ is the Whitehead link group, then $c = 1+i$.
\item If $\Gamma$ is the $6_2^2$ link group, then $c = \dfrac{3+i\sqrt{3}}{2}$.
\item If $\Gamma$ is the $6_3^2$ link group, then $c = \dfrac{1+i\sqrt{7}}{2}$.
\item If $\Gamma$ is a $\mathbb{Z}_2$-extension of the figure eight knot group, then $c = \sqrt{\dfrac{1+i\sqrt{3}}{2}}$.
\item If $\Gamma$ is a $\mathbb{Z}_2$-extension of the Whitehead link group, then $c = \sqrt{1+i}$.
\item If $\Gamma$ is a $\mathbb{Z}_2$-extension of the $6^2_2$ link group, then $c = \sqrt{\dfrac{3+i\sqrt{3}}{2}}$.
\item If $\Gamma$ is a $\mathbb{Z}_2$-extension of the $6^2_3$ link group, then $c = \sqrt{\dfrac{1+i\sqrt{7}}{2}}$.
\item If $[\mathrm{PSL}_2 (O_1) : \Gamma] = 8$, then $c = 1+i$.
\item If $\Gamma = \mathrm{PSL}_2 (O_3)$, then $c = \dfrac{1+i\sqrt{3}}{2}$.
\item If $[\mathrm{PSL}_2 (O_7) : \Gamma] = 2$, then $c = \dfrac{1+i\sqrt{7}}{2}$.
\item If $\Gamma = \mathrm{PGL}_2 (O_1)$, then $c = \dfrac{1+i}{\sqrt{2}}$.
\item If $[\mathrm{PSL}_2 (O_2) : \Gamma \cap \mathrm{PSL}_2 (O_2)] = 24$, then $c = \sqrt{2}+i$.
\item If $[\mathrm{PSL}_2 (O_3) : \Gamma \cap \mathrm{PSL}_2 (O_3)] = 30$, then $c = \dfrac{1+i\sqrt{3}}{\sqrt{2}}$.
\item If $\Gamma = \mathrm{PGL}_2 (O_3)$, then $c = \dfrac{\sqrt{3}+i}{2}$.
\item If $[\mathrm{PSL}_2 (O_{15}) : \Gamma \cap \mathrm{PSL}_2 (O_{15})] = 6$, then $c = \dfrac{\sqrt{3}+i\sqrt{5}}{2}$.
\end{itemize}
By Proposition \ref{arithcomp}, $J(\Gamma) = J(A,B) = |\tr [A,B] - 2| = |c|^2$; the result now follows.
\end{proof}

We conclude with a proposition that yields computations of J\o{}rgensen numbers of several two-bridge knot groups as a corollary.  For convenience, if $\mathbb{H}^3/\Gamma$ is the complement of a two-bridge knot $K$ in $S^3$, then the J\o{}rgensen number of $K$ is $J(K) = J(\Gamma)$.

\begin{pro} \label{knotcomp}
Let $\mathbb{H}^3/\Gamma$ be a hyperbolic 3-manifold of finite volume and $R^\ast$ the group of units in $R$, the ring of integers in $\mathbb{Q}(\tr\Gamma)$. Suppose $\tr \Gamma \subset R$ and \[R^\ast \cong W \times \langle u \rangle\] where $W$ is a finite cyclic group consisting of roots of unity and $u$ is a fundamental unit with $|u| > 1$. If $\Gamma = \langle A,B \rangle$ with $J(A,B) = |u|^2 < |\tr^2 X - 4|$ for all loxodromic elements $X$ in $\Gamma$, then $1 < \widetilde{J}(\Gamma) \leq J(\Gamma) = |u|^2$.
\end{pro}

\begin{proof}
Note that $\Gamma$ is not the figure-eight knot group since $R^\ast$ would then be finite. Hence, $1 < \widetilde{J}(\Gamma) \leq J(\Gamma)$ by Theorem \ref{part1}. Suppose $\Gamma = \langle X,Y \rangle$ with $J(X,Y) \leq |u|^2$. By hypothesis, $X$ must then be parabolic, so we have
\[1 < J(X,Y) = |\tr [X,Y] - 2| \leq |u|^2.\]
Similarly, $A$ must also be parabolic, so $J(A,B) = |\tr [A,B] - 2| = |u|^2$; thus, $\tr [A,B] - 2 = \xi_1 u^2$ for some root of unity $\xi_1$.

By Corollary \ref{unitmult}, $\tr [X,Y] - 2$ is a unit multiple of $\tr [A,B] - 2$ in $R$; hence, $\tr [X,Y] - 2 = \xi_2 u^a$ for some root of unity $\xi_2$ and integer $a$. But $|u| > 1$, so the inequality above forces $a = 1$ or 2.

Since $\tr X = \pm 2$, standard trace relations (see, for instance, 3.15 in Section 3.4 of \cite{MaclachlanReid03}) yield
\[ \tr [X,Y] - 2 = \tr^2 Y \pm 2\tr Y \tr XY + \tr^2 XY = (\tr Y \pm \tr XY)^2. \]
Thus, $\tr [X,Y] - 2$ is a square in $R^\ast$, so $a = 2$ and $J(X,Y) = |\tr [X,Y] - 2| = |u|^2$.  The result now follows.
\end{proof}

\begin{cor}
The J\o{}rgensen numbers of several two-bridge knots are as follows.
\[\begin{tabular}{c|c}
$ K $&$ J(K) $\\
\hline
$5_2 $&$ 1.32471796 $\\
$6_1 $&$ 1.55603019 $\\
$7_4 $&$ 2.20556943 $\\
$7_7 $&$ 1.55603019 $\\
\end{tabular}\]
\end{cor}

\begin{proof}
Let $K$ be one of the knots above, which are the two-bridge knots (7/3), (9/5), (15/11), and (21/13) respectively, and $\Gamma$ the Kleinian group with $\mathbb{H}^3/\Gamma = S^3 \smallsetminus K$.  Following the procedure outlined in Section \ref{bounds}, we see that Lemma \ref{tracefields} yields $\tr \Gamma \subset R$, the ring of integers in $\Q(\tr \Gamma) = \Q(z)$, where $z$ is as follows.

\[\begin{tabular}{c|c|c}
$ K $& minimum polynomial of $z$ & numerical value of $z$ \\
\hline
$5_2 $&$ 1+2z+z^2+z^3 $&$-0.21507985+1.307141279i$\\
$6_1 $&$ 1-2z+3z^2-z^3+z^4 $&$0.104876618-1.552491820i$\\
$7_4 $&$ 1+4z-4z^2+z^3 $&$2.10278472+0.665456952i$\\
$7_7 $&$ 1-z+3z^2-2z^3+z^4 $&$0.95668457-1.227185638i$\\
\end{tabular}\]

As verified by PARI (\cite{PARI2}), $J(A,W) = |z| = |u|^2$, where $A$ and $W$ are as in Section \ref{bounds} and $u$ is a fundamental unit in $R$. Thus, we see that $\mathbb{H}^3/\Gamma$ satisfies all hypotheses of Proposition \ref{knotcomp} but the last. To verify it, first suppose that $X$ is a loxodromic element in $\Gamma$ whose axis projects to a closed geodesic on $\mathbb{H}^3/\Gamma$ of length $L \geq 3$ (see, for instance, Section 1.3 of \cite{MaclachlanReid03}).  Since \[|\tr X| = 2\cosh\dfrac{L}{2}\] (Section 7.34 of \cite{Beardon83}) and hyperbolic cosine is an increasing function, $L \geq 3$ implies \[|\tr X| \geq 2\cosh\dfrac{3}{2} = 4.704819230.\]  Thus, $|\tr^2 X| \geq 22.13532399$, and so \[ |\tr^2 X - 4| \geq ||\tr^2 X| - 4| \geq 18.13532399.\]

Using SnapPea (\cite{SnapPea}, \cite{HodgsonWeeks94}), we now enumerate squares of traces of loxodromic elements in $\Gamma$ whose axes project to closed geodesics on $\mathbb{H}^3/\Gamma$ of length less than three to find the value of
\[ \alpha = \inf \{|\tr^2 X - 4| : X \in \Gamma \text{ is loxodromic} \}.\]  The values of $|z| = |u|^2$ and $\alpha$ for the aforementioned knots are as follows.

\[\begin{tabular}{c|c|c}
$ K $&$ |z| = |u|^2 $&$ \alpha $\\
\hline
$5_2 $&$ 1.32471796 $&$ 4.219276205 $\\
$6_1 $&$ 1.55603019 $&$ 3.955211258 $\\
$7_4 $&$ 2.20556943 $&$ 4.434378815 $\\
$7_7 $&$ 1.55603019 $&$ 5.105997169 $\\
\end{tabular}\]

Thus, in each case we have \[ |u|^2 < \alpha = \inf \{|\tr^2 X - 4| : X \in \Gamma \text{ is loxodromic} \}, \] so the result now follows from Proposition \ref{knotcomp}.
\end{proof}

\subsection*{Acknowledgements}
The author thanks his advisor, Alan Reid, whose patience and guidance made this work possible, and Grant Lakeland for helpful conversation regarding the fundamental polyhedra of Bianchi groups.

\bigskip

\address{
\noindent The University of Texas at Austin\\
Department of Mathematics\\
1 University Station C1200\\
Austin, TX 78712, USA}

\medskip

\noindent\email{callahan@math.utexas.edu}

\end{document}